\documentclass{amsart} 
\usepackage{amssymb}
\usepackage{url}
\usepackage{amscd}

\usepackage[mathscr]{euscript}

\newtheorem{thm}{Theorem}[section]
\newtheorem{cor}{Corollary}
\newtheorem{prop}[thm]{Proposition}
\newtheorem{lem}[thm]{Lemma}

\theoremstyle{remark}

\theoremstyle{definition}
\newtheorem{dfn}[thm]{Definition}
\newtheorem{ex}[thm]{Example}
\newtheorem{exs}[thm]{Examples}
\newtheorem{corollary}{Corollary}
\newtheorem{remark}[thm]{Remark}
\newtheorem{example}[thm]{Example}
\newtheorem{lemma}[thm]{Lemma}
\newtheorem{subth}{ }[thm]
\newtheorem{notation}[thm]{Notation}

\numberwithin{equation}{section}
\numberwithin{thm}{section}
\numberwithin{cor}{thm}
\numberwithin{sublem}{thm}
\numberwithin{rem}{thm}

\newcommand{\steq}{\stepcounter{equation}} 
\def\>{{\mkern 1mu}} 
\def\<{{\mkern-1mu}}

\newcommand{\set}{\!:=} 
\def\dash#1.#2.{$#1$-#2}

\newcommand{\Hom}{{\operatorname{Hom}}} 
\newcommand{\Ext}{\operatorname{Ext}} 
\newcommand{\Tor}{\operatorname{Tor}}

\newcommand{\Spec}{{\textup {Spec}}} 
\newcommand{\coker}{\operatorname{coker}}

\newcommand{\lra}{\longrightarrow} 
\newcommand{\xto}[1]{{\xrightarrow{#1}}} 
\newcommand{\iso}{{\mkern8mu\longrightarrow \mkern-25.5mu{}^\sim\mkern17mu}}

\newcommand{\Otimes}{\underset{\vbox to 0pt {\vskip-1.05ex\hbox{$\scriptscriptstyle=$}\vss}}\otimes}

\newcommand{\Rf}{{\R f_{\<*}}} 
\newcommand{\Lf}{{\Ll f^*}}

\newcommand{\tg}[1]{\tau^{}_{{\scriptscriptstyle \ge} #1}} 

\font\cmss=cmss8 
\font\cmsx=cmss8 scaled 720

\newcommand{\cmsc}{{\text {\cmss\char'143}}}

\newcommand{\cmsq}{{\text {\cmss\char'161}}} 
\newcommand{\qc}{{\cmsq\cmsc}} 
\newcommand{\qcc}{{\text{\cmsx\char'161 \cmsx\char'143}}}

\newcommand{\ScriptO}{{\mathcal O}} 
\newcommand{\SO}{{\mathcal O}} 
\newcommand{\D}{{\mathbf D}} 
\newcommand{\R}{{\mathbf R}} 
\newcommand{\Ll}{{\mathbf L}} 
\newcommand{\Z}{{\mathbb Z}}

\newcommand{\Dqc}{\mathbf D_{\qc}\mkern1mu}

\newcommand{\OX}{{\ScriptO_{\<\<X}}}

\def\smalldrlm#1{\underset{\vtop{\vskip-4.5pt\hbox to 1.37em{\rightarrowfill} \vskip-10pt\hbox{$\scriptscriptstyle \;\lift.6,{\ #1},$}}}\lim\,}

\def\lift#1,#2,{\vbox to 0pt{\vskip-#1 ex\hbox{$\scriptstyle #2$}\vss}} 
\def\vs#1{\vspace{#1pt}}

\newcommand\lemm[1]{\begin{lemma} \label{#1}} 
\newcommand\pro[1]{\begin{prop} \label{#1}} 
\newcommand\thmm[1]{\begin{thm} \label{#1}} 
\newcommand\corr[1]{\begin{corollary} \label{#1}}

\newcommand\sthm[1]{\begin{subth} \label{#1}} 
\newcommand\exm[1]{\begin{example} \label{#1} \begin{em}} 
\newcommand\plm[1]{\begin{prblm} \label{#1} \begin{em}} 
\newcommand\rmk[1]{\begin{remark} \label{#1}} 
\newcommand\rmd[1]{\begin{reminder} \label{#1} \begin{em}} 
\newcommand\ntn[1]{\begin{notation} \label{#1}}

\newcommand\elem{\end{lemma}} 
\newcommand\erfn{\end{refinement}} 
\newcommand\eprt{\end{proto}} 
\newcommand\ethm{\end{thm}} 
\newcommand\ecor{\end{corollary}} 
\newcommand\edfn{\end{dfn}} 
\newcommand\esthm{\end{subth}} 
\newcommand\epro{\end{prop}} 
\newcommand\etri{\end{triviality}} 
\newcommand\eexm{\end{em} \end{example}} 
\newcommand\efac{\end{em} \end{facts}} 
\newcommand\ecmt{\end{em} \end{com}} 
\newcommand\ermk{ \end{remark}} 
\newcommand\ermd{\end{em} \end{reminder}} 
\newcommand\eill{\end{em} \end{illustration}} 
\newcommand\eplm{\end{em} \end{prblm}} 
\newcommand\ecas{\end{em} \end{case}} 
\newcommand\ecau{\end{em} \end{caution}} 
\newcommand\ecax{\end{em} \end{cauex}} 
\newcommand\eimn{\end{em} \end{importnota}} 
\newcommand\entn{ \end{notation}} 
\newcommand\econ{\end{em} \end{construction}} 
\newcommand\esmr{\end{em} \end{summary}} 
\newcommand\ehyp{ \end{hypo}} 
\newcommand\ecnc{\end{em} \end{conclusion}} 
\newcommand\essthm{\end{em} \end{ssubth}}

\newcommand\script{\mathscr}

\newcommand\ct{{\script T}}

\newcommand\cs{{\script S}} 
\newcommand\cS{{\script S}}

\newcommand\ca{{\script A}}

\newcommand\bd{\begin{description}} 
\newcommand\ed{\end{description}}

\newcommand\prf{\begin{proof}} 
\newcommand\eprf{\end{proof}}

\newcommand\la{\longrightarrow}

\newcommand\eq{\quad=\quad} 
\newcommand\BA{\begin{array}{c}} 
\newcommand\EA{\end{array}}

\newcommand\mr\Modtc

\newcommand\ra\ri

\newcommand\Mod[1]{\ensuremath{\mathop{\textup{Mod-}#1}}\xspace} 
 
\newcommand\Modtc{\Mod{\ct^c}}

\newcommand\wt{\widetilde} 

\begin{document}

\date{Received November 24, 2006; received in final form March 16, 2007}

\title[Quasi-perfection and boundedness]{Quasi-perfect scheme-maps and
boundedness of the twisted inverse image functor}

\author{Joseph Lipman}
\address{Joseph Lipman\\
Department of Mathematics\\
Purdue University\\
W. Lafayette, IN 47907\\
USA}
\email{jlipman@purdue.edu}

\author{Amnon Neeman}
\address{Amnon Neeman\\
Centre for Mathematics and its Applications\\ 
Mathematical Sciences Institute\\
John Dednam Building\\
The Australian National University\\
Canberra, ACT 0200\\
Australia}
\email{Amnon.Neeman@anu.edu.au}

\dedicatory{To Phillip Griffith, on his 65th birthday}

\thanks{The first author is partially supported by  National Security
Agency award H98230-06-1-0010.  The second author is partially supported
by the Australian Research Council.}

\begin{abstract}
For a map $f\colon X\to Y$ of quasi-compact quasi-separated schemes, we discuss
 \emph{quasi-perfection}, i.e., the right adjoint~$f^\times\mkern-1mu$ of~ 
$\mkern1mu\mathbf Rf_*$ \mbox{respects} small direct sums.
This is equivalent to the existence of a functorial isomorphism 
$f^\times\mathcal O_{\mkern-2mu Y}\otimes^{\mathbf L} \mathbf Lf^*(\<-\<)\!
{\mkern8mu\longrightarrow \mkern-25.5mu{}^\sim\mkern10mu} f^\times (-)$;  
to \emph{quasi-properness}
(preservation by~$\Rf$ of pseudo-coherence, or just \emph{properness} in the noetherian\- case) plus  boundedness of  $\mathbf Lf^*\<$ (finite tor-dimensionality), 
or of the functor~$f^\times\<$; and to some other conditions. We use a globalization, previously known only for divisorial schemes, of
the local definition of pseudo-coherence of complexes,
as well as a refinement of the known fact that the derived 
category of complexes with quasi-coherent homology is generated by a single
perfect complex.
\end{abstract}


\maketitle


\section{Introduction}
\label{Intro}
This paper, inspired by \cite[p.\,396, Lemma 1 and Corollary 2]{V},  deals with matters	 raised there, but not yet fully treated in the literature.

Throughout, \emph{scheme} will  mean \emph{quasi-compact
quasi-separated scheme} (see \cite[\S6.1, p.\,290ff]{GD}),
though weaker assumptions would sometimes
suffice. Unless\- otherwise indicated, a~\emph{map} $f\colon X\to Y$ will be a \emph{scheme\kern.5pt-morphism,}
necessarily quasi-compact and quasi-separated.

For a scheme $X\<\mkern-.5mu$, $\D(X)$ is the (unbounded) derived category
of the category of (sheaves of) \dash\OX.modules., and $\Dqc(X)$  is
the full subcategory
whose objects\- are the  \dash\OX.complexes. whose homology sheaves are 
all quasi-coherent. For any map $f\colon X\to Y$, 
 the
derived functor \mbox{$\Rf\colon\D(X)\to\D(Y)$} takes $\Dqc(X)$ to~$\Dqc(Y)$ \cite[Prop.\,(3.9.2)]{Lp}. 
Grothendieck Duality theory\vspace{.7pt} asserts, to begin,  that \emph{the restriction\/ 
\mbox{$\Rf\colon\Dqc(X)\to \Dqc(Y)$}\vspace{1.2pt} has a right~adjoint~
$f^\times\<$}, the ``twisted inverse image functor" in our title.\vs{-1}%
\footnote{
Warning: for nonproper maps of noetherian schemes, 
the usual twisted inverse image~$f^!$ 
differs from $f^\times\<$, and is not covered by this paper. For that, see, e.g., \cite[\S4.9]{Lp}.}

A~proof for maps of  \emph{separated} schemes, 
suggested by Deligne's appendix to~\cite{H}, is described in \cite[\S4.1]{Lp}. This proof
depends ultimately on the Special Adjoint\- Functor Theorem, 
applied to categories
of sheaves. A more direct approach, via Brown Representability---which
applies immediately to derived categories---is given in \cite{N1}.
Originally this too required separability, but now that assumption can
be dropped because of \cite[p.\,9, Thm.\,3.3.1]{BB}, which gives that 
$\Dqc(X)$ is compactly generated, and because $\Rf$ commutes with
$\Dqc$-coproducts (=\;direct sums)\vs1 \cite[(3.9.3.3)]{Lp}.%
\footnote{Subsequently, a slightly simpler proof was given in \cite[p.\,14, 3.3.4]{BB}.
(In that proof one needs to replace ``flabby'' by ``quasi-flabby,''
see \cite[\S2]{Kf}.)
}

The functor $f^\times$ emerging from these proofs commutes with  translation (=suspension) of complexes, and is \emph{bounded-below}
(\emph{way-out right} in the sense of \mbox{\cite[p.\,68]{H}),}
i.e., there exists 
an integer $m$
such that for~every $F\in\Dqc(Y)$ with $H^iF=0$ for all $i$ 
less than some 
integer~$n(F)$, it holds that
$H^i\<f^\times\<\< F=0$ for all $i<n(F)-m$ (see \cite[(4.1.8) and
the remarks preceding it]{Lp}).\looseness=-1

``Bounded-below" has a similar
meaning for any functor between derived categories.
\emph{Bounded-above} is defined in an analogous way, with $>$
(resp.~$+$) in place of $<$ (resp.~$-$). 
A functor is \emph{bounded} if it is bounded both above and below. 
Boundedness  enables a potent form of induction in derived categories, expressed by the ``way-out Lemmas"  \cite[p.\,68, Prop.\,7.1 and~p.\,73, Prop.\,7.3]{H}.\looseness=-1

For example, the left adjoint $\Lf\<$ of $\>\Rf$ is always bounded-above; and $\Lf$ is bounded iff $f$ has \emph{finite tor-dimension} (a.k.a. \emph{finite flat dimension}), 
that is, there is an integer $d\ge0$ such that for each $x\in X$ there exists an
exact sequence of \dash\ScriptO_{Y,\>f(x)}.modules.
$$
0\to P^{}_d\to P^{}_{d-1}\to\dots\to P^{}_{\>1}\to P^{}_{\>0}\to\ScriptO_{\<\<X\<\<,\>x}\to 0
$$
with $P^{}_i$ flat over $\ScriptO_{Y,\>f(x)}$ $(0\le i\le d\>).$\vspace{2pt}

We will be concerned with the relation between boundedness of 
the right adjoint~$f^\times$ and the left adjoint $\Lf\<\<$,
\kern.5pt especially in the context of \emph{quasi-perfection,} a property of maps to be discussed at length now and in \S2.%
\vspace{1pt}\looseness=-1

\begin{dfn}\label{quasiperfect}
We say a map $f\colon X\to Y$ is \emph{quasi-perfect} if $f^\times$~ respects direct sums  in~
$\Dqc$, i.e., for any small
\dash\Dqc(Y).family.~$(E_\alpha)$
the natural map is an isomorphism
$$
\underset{\lift.9,\alpha,}\oplus\>\> f^\times\<\< E_\alpha\iso f^\times\<
(\underset{\lift.9,\alpha,}\oplus\>\> E_\alpha).
$$
\end{dfn}
As will be explained below, quasi-perfection is also characterized by the existence of a canonical isomorphism
$$
\smash{f^\times\SO_Y\mkern-.5mu\otimes^{\mathbf L} \Lf\<\< F\iso f^\times\<\< F}\qquad\bigl(F\in\Dqc(Y)\bigr).\vspace{-1pt}
$$

More characterizations are given in \S2---for instance, via compatibility of~$f^\times$ with tor-independent base change (Theorem~\ref{Thm2.7}). That  section also brings in the
related condition on maps of being \emph{perfect}, i.e.,
pseudo\kern.5pt-coherent and of finite
tor-dimension. (Pseudo-coherence will be reviewed in~\S2. It holds, for instance, for all finite-type maps of noetherian schemes; and then descent to the noetherian case yields  that
every flat, finitely presentable map~is \mbox{pseudo\kern.5pt-coherent.)} For example, for a proper map\/~$f$ of noetherian schemes, \emph{$f$~is~quasi-perfect\/ 
$\Leftrightarrow f$ is perfect\/ $\Leftrightarrow f^\times\!$ is bounded.}\vs1

It is stated in \cite[p.\,396, Lemma 1]{V} that any proper map $f$ of finite-dimensional noetherian schemes is quasi-perfect. In general, however, this fails even for closed immersions.
But $f^\times$ \emph{does} respect direct sums when the summands $E_\alpha$ are uniformly homologically
bounded below, i.e., there exists an integer~$n$ such that for all $\alpha$, 
$H^i\<E_\alpha=0$ whenever $i< n$ \cite[(4.7.6)(b)]{Lp}. Consequently, 
\emph{if~the functor $f^\times\!$ is bounded, then $f$ is quasi-perfect.}

Our main results say more.  But first, call a map 
$f\colon X\to Y$ \emph{quasi-proper} 
if $\Rf$ takes pseudo\kern.5pt-coherent $\OX$-complexes to 
pseudo\kern.5pt-coherent $\SO_Y$-complexes.
(Again, pseudo-coherence is explained in  \S2. In particular, if $X$ is noetherian then
$E\in\mathbf \D(X)$ is pseudo\kern.5pt-coherent iff the homology sheaves~$H^n(E)$ are all coherent, and vanish for $n\gg0$.)  Kiehl showed that \emph{every proper pseudo-coherent map is quasi-proper.} Consequently, 
any flat,  finitely presentable, proper map, being pseudocoherent, is
quasi-proper (and perfect and quasi-perfect as well). Moreover, 
\emph{when\/ $Y\<$ is noetherian, every finite-type separated 
quasi-proper\/ \mbox{$\>f\colon X\to Y$} is proper.}\vspace{2pt}

Here are the main results.

\begin{thm}\label{Thm1.2} 
For a map\/ $f\colon X\to Y,$ the following are equivalent\kern1pt$:$

{\rm(i)} $f$ is quasi-perfect \textup{(}resp.~perfect\/\textup{)}.

{\rm(ii)} $f$ is quasi-proper \textup{(}resp.~pseudo-coherent\/\textup{)} and has finite tor-dimension.

{\rm(iii)} $f$ is quasi-proper \textup{(}resp.~pseudo-coherent\/\textup{)} and  $f^\times\!$ is bounded.
\end{thm}

Hence, by Kiehl's theorem,  \emph{every proper perfect map is quasi-perfect}.\vs2

The implication (i)${}\Rightarrow{}$(iii) is worked out in \S4. The proofs in \S4 are based on Theorems~\ref{perfapprox} and ~\ref{perfgen}, which
are of independent interest.

Theorem~\ref{perfapprox} states that for a scheme~$X\<$, any 
pseudo\kern.5pt-coherent $\OX$-complex can be 
``arbitrarily-well approximated,"
\emph{globally,} by a perfect complex. (\emph{Local} approximability is essentially the definition
of pseudo-coherence. The global result was previously known only for divisorial schemes.\kern-1pt) 

This leads to quasi-proper maps being characterized as those $f$ such that $\Rf$ takes perfect complexes to pseudo\kern.5pt-coherent ones. 
Since by Prop.\,\ref{Prop2.1}, quasi-perfect maps are those $f$ such that
$\Rf$ takes perfect complexes to perfect ones, it follows at once that 
quasi-perfect maps are quasi-proper.\vs1

Theorem~\ref{perfgen} refines a theorem of Bondal and 
van den Bergh \cite[p.\,9, Thm.\,3.1.1]{BB}
which states that
the triangulated category $\Dqc(X)$ is generated by a single perfect complex. 
With this in hand, one can prove  Corollary~\ref{qp bounded}, which says 
that for any 
quasi-perfect \emph{or} perfect  $f$ as above,
$f^\times\!$ is bounded.

The implication (iii)${}\Rightarrow{}$(ii) results from
Theorem~\ref{Thm3.1}, which says, for any $f\colon X\to Y$ as above, 
\emph{if\/ $f^\times\!$ is bounded then\/ $f$ has finite
tor-dimension.}

Finally, the implication (ii)${}\Rightarrow{}$(i) holds by definition for the resp.~case, and  is proved for the other case in \S2, Example~\ref{Example}(a).\vs4

Let us call a map $f\colon X\to Y$  \emph{locally embeddable} if  every $y\in Y$ has an open neighborhood $V$ over which  the\vspace{.5pt}
induced map  \mbox{$f^{-1}V\to V$} factors as
$f^{-1}V\,\xto{\lift.8,i,}\,Z\,\xto{\lift.65,p,}\,V$ where $i$ is  a closed
immersion and $p$~is~smooth. (For~instance, any~quasi-projective  $\<f\<$ satisfies this condition.)
Proposition~\ref{Prop2.5} asserts that \emph{any quasi-proper locally
embeddable map   
is pseudo-coherent.} A similar proof shows that any quasi-perfect locally embeddable map is perfect. By 1.2, then,  \emph{a locally embeddable map
is quasi-perfect iff it~is quasi-proper and~perfect.}\vs3

The equivalence of (i) and (ii) in Theorem~\ref{Thm1.2}
generalizes  \cite[p.\,396, Cor.\,2]{V}, in view of the following 
characterization
(mentioned above) of quasi-perfection.\vs1

For a map $f\colon X\to Y\<$, and for
any \mbox{$E\in\Dqc(X)$,}
$F\in\Dqc(Y)$, with \smash{$\Otimes\set\otimes^{\mathbf L},$} the derived tensor product, the ``projection map"\vspace{-2pt}
$$
\pi\colon\smash{(\Rf E)\Otimes F \to \Rf(E\Otimes\Lf F),}
$$
defined to be adjoint to the natural composite map
$$
\smash{\Lf\bigl((\Rf E)\Otimes F\bigr)\iso(\Lf\Rf E)\Otimes\Lf F\to E\Otimes \Lf F,}
$$
is an \emph{isomorphism.}
(This is well-known under more restrictive hypotheses; for a proof in the stated generality, see \cite[Prop.\,(3.9.4)]{Lp}.) 
There results a natural map\vs{-2}
\steq\steq
\begin{equation}\label{chi}
\smash{\chi^{}_F\colon f^\times\SO_Y\mkern-.5mu\Otimes \Lf\<\< F\to f^\times\<\< F}\qquad\bigl(F\in\Dqc(Y)\bigr),
\end{equation}
adjoint to the natural composite map
$$
\smash{\R f_* (f^\times\< \SO_Y\<\<\Otimes \Lf \<\<F)}
\underset{\lift1.2,\<\<\pi^{-\<1},}{\iso}
\R f_*  f^\times\smash{\< \SO_Y\<\<\Otimes F}
\to \SO_Y\mkern-.5mu\Otimes F =F.
$$
It is clear (since \smash{$\Otimes$} and  $\Lf$ both respect direct sums, see e.g., \cite[3.8.2]{Lp}) that \emph{if\/ $\chi^{}_F$ is an isomorphism for all\/ $F\in\Dqc(Y)$ then $f$ is quasi-perfect;} and Proposition~\ref{Prop2.1} gives the converse.\vs{-1}

\section{quasi\kern.5pt-perfect maps}
\label{section1}

For surveying quasi-perfection in more detail, starting with Proposition~\ref{Prop2.1},
we need some preliminaries. 

First, a brief review of the notion of pseudo-coherence of complexes.
(Details can be found in the primary source~\cite[Expos\'e III]{I}, or, perhaps more accessibly, in \cite[pp.\;283{\it ff,} \S2]{TT}; a summary appears in 
\mbox{\cite[\S4.3]{Lp}.)} The idea is built up from that of
\emph{strictly perfect} $\OX$-complex, i.e., bounded complex 
of finite-rank free \dash\OX.modules..

For $n\in\mathbb Z$, a map $\xi\colon P\to E$ in~ 
$\mathbf K(X)$, the homotopy category 
of $\OX$-complexes, (resp.~in $\mathbf D(X)$), is said to be an
$n$-\emph{quasi-isomorphism} (resp. \mbox{$n$-\emph{isomorphism}}) if the
following two equivalent conditions hold:\vs1

(1) The homology map $H^j(\xi)\colon H^j(P)\to H^j(E)$ is bijective 
for all $j>n$ and surjective for $j=n$.\vs{.6}

(2) For any $\mathbf K(X)$- (resp.~$\D(X)$-)triangle
$$
\CD
P@>\xi>>E @>>> Q @>>> P[1],
\endCD
$$
it holds that $H^j(Q)=0$ for all $j\ge n$.\vs1

Then $E$ is said to be $n$-\emph{pseudo-coherent} if $X$ has an open
covering $(U_\alpha)$ such that for each $\alpha$ there exists a
strictly perfect $\SO_{U_\alpha}$-complex~$P_\alpha$ and an 
$n$-quasi-isomorphism (or equivalently, an $n$-isomorphism) 
$P_\alpha\to E|_{U_\alpha}$, see \cite[p.\,98, D\'efinition 2.3]{I};
and $E$ is \emph{pseudo-coherent} if $E$ is $n$-pseudo\kern.5pt-coherent for every
$n$.\looseness=-1

 \noindent If $\OX$
is coherent, this means simply that $F$~has coherent homology sheaves,
vanishing in all sufficiently large degrees
 [{\it ibid.,} p.\,116, top].
When $X$ is noetherian and finite-dimensional, it means
that $F$~is globally \dash\D.isomorphic. to a bounded-above
complex of coherent \dash\OX.modules. [{\it ibid.,} p.\,168, 
Cor.\,2.2.2.1].

\smallskip

A complex $E\in\mathbf D(X)$ ($X$ a scheme) is said to be \emph{perfect} if
it is locally
\mbox{\dash\mathbf D.isomorphic.} to a strictly perfect $\OX$-complex. 
More precisely, $E$ is said to have \emph{perfect amplitude in} $[a,b\>]$
($a\le b\in\mathbb Z$) if locally on~$X\<$, $E$~is \dash\mathbf
D.isomorphic. to a bounded complex of finite-rank  free
 \dash\OX.modules.
vanishing in all degrees $<a\>$ or $>b$. Thus $E$ is perfect
iff it has perfect amplitude in some interval~$[a,b\>]$.\looseness=-1 

By \cite[p.\,134, 5.8]{I}, $E$ has perfect amplitude in $[a,b\>]$  iff $E$ is 
$(a-1)$-pseudo\kern.5pt-coherent and  has tor-amplitude in $[a,b\>]$ (i.e., is
globally \dash\mathbf D.isomorphic. to a flat complex
vanishing in all degrees \hbox{$<a$ and $>b$).} So 
$E$ \emph{is perfect iff it is pseudo-coherent and has finite tor-dimension}
(the latter meaning that it is \dash\mathbf D.isomorphic. to a bounded flat complex).

\smallskip

A  map $f\colon X\to Y$ is \emph{pseudo\kern.5pt-coherent} if 
every  $x\in X$ has an open neighborhood~$U$ such that the restriction $f|_U$ factors as
$\>U\>\xto{\lift.8,i,}\>Z\>\xto{\lift.65,p,}\>Y\<$, where  $i$ is a closed
immersion such that $i_*\SO_U$ is pseudo\kern.5pt-coherent on $Z$, and
$p$ is  smooth \cite[p.\,228, D\'efn.\,1.2]{I}. Pseudo-coherent
maps are finitely presentable. Compositions of
pseudo\kern.5pt-coherent
maps are pseudo\kern.5pt-coherent \cite[p.\,236, Cor.\,1.14]{I}.

\smallskip
A map is \emph{perfect} if it is pseudo\kern.5pt-coherent and has finite tor-dimension
\cite[p.\,250, D\'efn.\,4.1]{I}. Any smooth map is perfect, any regular immersion
(=\;closed immersion corresponding to a quasi-coherent ideal generated locally\-
by a regular sequence) is perfect, and compositions of perfect
maps are perfect \cite[p.\,253, Cor.\,4.5.1(a)]{I}.\looseness=-1

For noetherian $Y\<$, any finite-type $f\colon X\to Y$ is pseudo\kern.5pt-coherent.
Pseudo\kern.5pt-coherence (resp.~ perfection) of maps survives  tor-independent base change \cite[p.\,233, Cor.\,1.10; p.\,257, Cor.\,4.7.2]{I}. Hence, by descent to the noetherian
case~\cite[IV, (11.2.7)]{EGA}, \emph{every flat 
finitely-presentable map is perfect.}\vs1

\smallskip

A map $f\colon X\to Y$ is \emph{quasi-proper} 
if $\>\Rf$ takes
pseudo\kern.5pt-coherent \dash\OX.complexes. to pseudo\kern.5pt-coherent
\dash\SO_Y.complexes..
\vspace{.7pt}

Kiehl's Finiteness Theorem \cite[p.\,315, {Thm.\,$2.9'$}]{Kl}
(first proved by Illusie
for projective maps \cite[p.\,236, Thm.\,2.2]{I}) generalizes preservation of coherence by higher direct images under
proper maps of noetherian schemes. It states that
\emph{every proper pseudo-coherent map is  quasi-proper.}

This theorem (or its special case \cite[p.\,240, Cor.\,2.5]{I}), plus \cite[Ex.\,(4.3.9)]{Lp}) implies that \emph{if\/ $Y\<$ is \mbox{noetherian} then a finite-type separated $f\colon X\to Y$ is  quasi-proper iff it is proper.}

\smallskip
For details in the  proof of the following Proposition, and for some subsequent considerations, recall that an object $C$ in a triangulated category~$\ct$ is \emph{compact\-} \kern.7pt if for every small $\ct$-family ~$(E_\alpha)$ the natural~map is an isomorphism\vs{-1}
$$
\underset{\lift.9,\alpha,}\oplus\>\> \Hom(C,E_\alpha)\!\iso\! 
\Hom(C,\underset{\lift.9,\alpha,}\oplus\>\> E_\alpha).\vs{-1.5}
$$

 For any scheme $X\<$,
the compact objects of $\Dqc(X)$ are just the perfect complexes, of which one is a generator \cite[p.\,9, Thm.\,3.1.1]{BB}.

\begin{prop}\label{Prop2.1}
For a map\/
$f\colon X\to Y\<,$ the following\vs1  are equivalent\/\kern.5pt$:$

{\rm(i)} $f$ is quasi-perfect \textup{(}Definition~\textup{\ref{quasiperfect}}\textup{)}.

{\rm(ii)} The functor $\R f_*$\ takes perfect complexes to perfect
complexes.

{\rm(ii)$'$} If\/ $S$ is a perfect generator of $\Dqc(X)$ then\/ $\Rf S$ is perfect.

{\rm(iii)} The twisted inverse image functor $f^\times$ has a right adjoint.

{\rm(iv)} For all\/ $F\in\Dqc(Y),$ the  map\/ in \eqref{chi} is an \emph{isomorphism}\vs{-2}
$$
\chi^{}_F\colon f^\times \SO_Y\<\<\Otimes \Lf\! F\iso f^\times\<\<F. 
$$

\end{prop}
\begin{proof} (i)${}\Leftrightarrow{}$(ii)${}\Leftrightarrow{}$(ii)$'$: \cite[p.\,224, Thm.\,5.1]{N1}.

(i)${}\Rightarrow{}$(iii): \cite[p.\,223, Thm.\,4.1]{N1}.

(iii)${}\Rightarrow{}$(i): simple.

(i)${}\Rightarrow{}$(iv)${}\Rightarrow{}$(i):  See \cite[p.\,226, Thm.\,5.4]{N1}.\vs1

To be precise, the results in \cite{N1} are proved for \emph{separated} schemes;
but with the remark preceding Prop.\,2.1, one readily verifies that the proofs survive without any separability requirement. 
\end{proof}

\begin{exs}\label{Example}   (a)  
 \emph{Any quasi-proper map\/ $f$ of finite tor-dimension---in particular, by Kiehl's theorem, any proper perfect map---is quasi-perfect.}  Indeed, $\Rf$ preserves
pseudo\kern.5pt-coherence, and by  \cite[p.\,250, 3.7.2]{I}
(a~consequence of the projection isomorphism mentioned near the end of the above Introduction), $\Rf$ preserves finite tor-dimensionality of
complexes; so  Prop.~\ref{Prop2.1}(ii) holds.\vs2\looseness=-1

(b)  Let $f\colon X\to Y$ be a map with $X$ 
\emph{divisorial}---i.e., $X$ has
an ample family $(\mathcal L_i)_{i\in I}$ of 
invertible \dash\OX.modules.  
\cite[p.\,171, D\'efn.\,2.2.5]{I}. Then \cite[p.\,212, Example 1.11 and p.\,224, Theorem 5.1]{N1} show that\vspace{.5pt} 
 $f$ \emph{is quasi\kern.5pt-perfect${}\Leftrightarrow{}$for each\/ $i\in I,$ the
\mbox{\dash\SO_Y.complex\/. $\Rf (\mathcal L_i^{\otimes -n_i})$} is perfect for all}\vs2
$n_i\gg 0$.\looseness=-1

(c)  Let $f\>$ be quasi-projective  and let $\mathcal L$ be\vspace{.5pt} 
 an \dash f.ample. invertible sheaf. Then
$f$ \emph{is quasi-perfect ${}\Leftrightarrow{}$ the\/
\dash\SO_Y.complex\/. $\R f_*(\mathcal L^{\otimes -n})$ is perfect for all}
$n\gg 0$. \vspace{.7pt}

Indeed, condition (ii) in Prop.\,\ref{Prop2.1}, together with the compatibility
of $\Rf$ and open base change, implies that quasi-perfection is a property of $f$ which 
can be checked locally on $Y\<$, and the same holds for perfection of $\R f_*(\mathcal L^{\otimes -n})$; so we may assume $Y$ affine, and apply (b).

\pagebreak[3]

(d) For a \emph{finite} map $f\colon X\to Y$ the following are equivalent:
\vspace{1pt}

(i) \kern2pt $f$ is quasi\kern.5pt-perfect.

(ii) $f$ is perfect.

(iii) The complex~$f_*\OX\cong\Rf\OX$ is perfect.\vspace{1.5pt}

This follows quickly from (a) and from Proposition~\ref{Prop2.1}(ii).
\end{exs}

\steq\steq
A \emph{tor-independent square} is a fiber square of maps
\begin{equation}\label{indt}
\CD
X'@>v>>X  \\
@V g VV  @VV f V  \\
Y' @>>\lift1.2,u,> Y\\[-3pt]
\endCD
\end{equation}
(that is, the natural map is an isomorphism $X'\iso X\times_Y Y'$) such that for all
$x\in X,$  $y'\in Y'$ and $y\in Y$ with $f(x)=u(y')=y$, and all $i>0$,
$\Tor_i^{\SO_{Y\!,y\>}}\!(\SO_{\<\<X\<,\>x}\>, 
\SO_{Y^{\lift.8,\scriptscriptstyle\<\prime,}\!,\>y'})=0$.

\medskip
The following stability properties will be useful. 
\stepcounter{thm}
\begin{prop}\label{Prop2.4} For any tor-independent 
square\/ \eqref{indt}$,$

{\rm(i)} If\/ the functor  $f^\times\!$ is bounded  then so is\/
$g^\times\<.$

{\rm(ii)} If\/ $f$ is quasi-perfect then so is\/ $g$.

{\rm (iii)} If\/ $f$ is quasi-proper then so is\/ $g$.
\end{prop}

\prf (i) and (ii) are proved in \cite[(4.7.3.1)]{Lp}; and (iii) is treated in Prop.\,\ref{stableqp} below (a slight change in whose proof gives another proof of (ii)).
\eprf

Since perfection (resp.~pseudo-coherence) is a local property of complexes, and $\Rf$ is compatible with open base change on $Y\<$, we deduce: 

\begin{cor}\label{Cor2.4.1} Let\/ $f\colon X\to Y$ be a map, and let\/ $(Y_i)_{i\in I}$ be an open cover of\/ $Y\<$. Then\/ $f$ is quasi-perfect \textup{(}resp.~quasi-proper\/\textup{)} $\Leftrightarrow$ for all\/ $i,$ the same is true of the induced map\/ $f^{-1}Y_i\to Y_i$.
\end{cor}

\begin{prop}\label{Prop2.5} Let\/ $f\colon X\to Y$ be a \emph{locally embeddable} map, i.e.,  every\/ $y\in Y$ has an open neighborhood\/ $V$ over which  the\vspace{.5pt}
induced map\/  \mbox{$f^{-1}V\to V$} factors as\/
$f^{-1}V\,\xto{\lift.8,i,}\,Z\,\xto{\lift.65,p,}\,V$ where $i$ is  a closed
immersion and $p$~is~smooth. \textup{(}\kern-.5pt For~instance, any~quasi-projective  $\<f\<$ satisfies this condition \textup{\cite[II, \kern-.5pt (5.3.3)]{EGA}}.\textup{)}\vs1

{\rm(i)} If $f$ is quasi-proper then $f$ is pseudo-coherent.

{\rm(ii)} If $f$ is quasi-perfect then $f$ is perfect.

\end{prop}

\begin{proof}

By Corollary~\ref{Cor2.4.1}, quasi-properness (resp.~quasi-perfection) of $f$ is a property local over $Y\<$; and since they are compatible with tor-independent base change, the same is true of pseudo-coherence and perfection. So we may as well 
assume that $X=f^{-1}V$.
Then  it  suffices to show that the complex $i_*\OX$ is pseudo\kern.5pt-coherent when $f$~is quasi-proper,\vs{.5} (resp., by \cite[p.\,252, Prop.\,4.4]{I}, 
that $i_*\OX$ is perfect when $f$~is quasi-perfect). 

But 
$i$ factors as $X\,\xto{\lift.6,\gamma,}\, X\!\times_Y\!Z\,\xto{\lift.6,g,} \,Z$ 
with $\gamma$  the graph of~$\>i$ and $g$ the projection. The map~$\gamma$ is
a local complete intersection \cite[IV, (17.12.3)]{EGA}, so
the complex $\gamma_*\OX$ is perfect. 
Also, $g$~arises from $f$ by flat base change, so by Proposition~\ref{Prop2.4}, $g$ is quasi-proper (resp.~quasi-perfect). Hence $i_*\OX\!=\R i_*\OX\!=\R g_*\<\gamma_*\OX\<$
is indeed pseudo\kern.5pt-coherent (resp.~perfect).
\end{proof}

\steq\steq
{\bf(2.6).} For any tor-independent square\/ \eqref{indt}, the map
\begin{equation*}\label{2.6.1}
\theta(E)\colon\Ll u^*\Rf E\to \R g_*\Ll v^*\<\<E \qquad\bigl(E\in\Dqc(X)\bigr)
\tag{2.6.1}
\end{equation*}
adjoint to the natural composition
$$
\Rf E\to \Rf\R v_*\Ll v^*\<\<E\cong \R u_*\R g_*\Ll v^*\<\<E
$$
(equivalently, to $\Ll g^*\Ll u^*\Rf E\cong\Ll v^*\Lf\Rf E\to\Ll v^*\<\<E$)
is an isomorphism, so that one has a \emph{base-change map}

\begin{equation*}\label{2.6.2}
\beta(F)\colon \Ll v^*\<\<f^\times\<\<F\to g^\times\Ll u^*\<\<F\qquad \bigl(F\in\Dqc(Y)\bigr)
\tag{2.6.2}
\end{equation*}
adjoint to the natural composition
$$
\R g_*\Ll v^*\<\<f^\times\<\<F\underset{\theta^{-1}}{\iso}\Ll u^*\Rf f^\times\<\<F\to \Ll u^*\<\<F\>.
$$

The fundamental \emph{independent base-change} theorem states that:\vspace{1pt} 

\emph{Let there be given a tor-independent square\/~\eqref{indt} and an\/
$F\in\Dqc(Y)$. If $f$ is  quasi-proper,
$u$ has finite tor-dimension,} 
\emph{and  $H^n\<F=0$ for all\/~$n\ll 0,$ then\/
$\beta(F)$~is an isomorphism.}\vspace{1.5pt} 

This theorem is well-known, at least under more restrictive hypotheses. For a treatment
in full generality, see \cite[\S\S4.4--4.6]{Lp}. \vspace{1pt}

One consequence, in view of Proposition~\ref{Prop2.4}(i), is:

\stepcounter{thm}\stepcounter{cor}\stepcounter{cor}
\begin{cor}\label{Cor2.6.3} Let\/ $f\colon X\to Y$ be a quasi-proper
map and let\/ $(Y_i)_{i\in I}$ be an open cover of\/ $Y\<$. Then\/ $f^\times\!$ is bounded  $\Leftrightarrow$ for all\/~$i,$ the same is true of the induced map\/ $f^{-1}Y_i\to Y_i\>$.
\end{cor}

For quasi-perfect $f\<$, a  stronger base-change theorem holds---which, together with
boundedness of $f^\times$ (Corollary~\ref{qp bounded}), characterizes quasi-perfection:

\begin{thm}[{\cite[Thm.\,4.7.4]{Lp}}]\label{Thm2.7}
Let
$$
\begin{CD}
X'@>v>>X  \\
@V g VV  @VV f V  \\
Y' @>>\lift1.2,u,> Y
\end{CD}
$$
be a tor-independent square, with $f$ 
quasi\kern.5pt-perfect.
Then for all $F\in\Dqc(Y)$ the base-change map  of\/~\textup{(2.6.2)} is an\/
{\rm isomorphism}\vs{-1}
$$
\beta(F)\colon  v^*\<\<f^\times\<\<F\iso g^\times u^*\<\<F.
$$
The same holds, with no assumption on $f,$\vadjust{\kern1pt} whenever\/ $u$ is finite and perfect.

\pagebreak[3]

Conversely, the following conditions on a map $f\colon X\to Y$ are equivalent;
and if  $f^\times\<$\ is bounded above, 
they imply that $f$~is quasi-perfect\/$:$\vs1

{\rm(i)} For any flat affine universally bicontinuous map  $u\colon Y'\to Y,\mkern.5mu$%
\footnote{\emph{Universally bicontinuous} means that for any $Y''\to Y$ the resulting projection map \mbox{$Y'\times_Y Y''\to Y''$} is a homeomorphism onto its image \cite[p.\,249, D\'efn.\,(3.8.1)]{GD}.}
the base-change map associated to the \textup{(}tor-independent\/\textup{)} square
$$
\CD
Y'\times_Y X={}@.X'@>v>>X  \\
@.@V g VV  @VV f V  \\
@.Y' @>>\lift1.2,u,> Y\vadjust{\kern-1pt}
\endCD
$$
is an isomorphism
$$
\beta(\SO_Y)\colon v^*\<\<f^\times\SO_Y\iso g^\times u^*\SO_Y.\vadjust{\kern1pt}
$$

{\rm (ii)} The map  in \eqref{chi} is an isomorphism
$$
\chi^{}_F\colon f^\times\SO_Y\Otimes \Lf\<\< F\iso f^\times\<\< F
\vadjust{\kern-4pt}
$$
 whenever $F$ is a flat quasi-coherent\/ \dash\SO_Y.{\kern.5pt}module..
\end{thm}

Keeping in mind Corollary~\ref{qp bounded} below ($f$ quasi-perfect${}\Rightarrow f^\times\<$ bounded), we can deduce:

\begin{cor}\label{2.7.1} A map\/  $f\colon X\to Y$  is 
quasi\kern.5pt-perfect\/  $\Leftrightarrow f^\times\!$\vspace{1pt} is 
bounded and the
following two conditions hold\/$:$\vs{1.5}

\textup{(i)} If\/ $u\colon Y'\to Y$ is an open immersion, and if\/ $v\colon Y'\times_Y X\to X$ and
$g\colon Y'\times_Y X\to Y$ are the projection maps, then the base-change map is an isomorphism\-\vspace{-2pt}
$$
\beta(\SO_Y)\colon v^*\<\<f^{\<\times}\<\SO_Y\iso g^\times u^*\SO_Y.
$$ 
Equivalently \textup{(}see \textup{[Lp, \S4.6, subsection V\kern.5pt]),}
for all $E\in\Dqc(X)$ the natural composite map 
$$
\Rf\mathcal \R\mathcal Hom^\bullet_{X}\<(E,\<f^{\<\times}\<\SO_Y\<\kern-.8pt)\<\<\to\<\<
\R\mathcal Hom^\bullet_{Y}\<(\Rf E,\Rf f^{\<\times}\<\SO_Y\<\kern-.8pt)\<\<\to
 \R\mathcal Hom^\bullet_{Y}\<(\Rf E,\SO_Y\<\kern-.8pt)
$$
is an isomorphism.\vspace{1.5pt}

\textup{(ii)} If\/ $(F_\alpha)$ is a filtered direct system of flat quasi-coherent\/
$\SO_Y$-modules,
then for all\/~$n\in\Z$ the natural map is an isomorphism
$$
\smalldrlm{\lift 0,\!\alpha{\,},} H^n(f^\times \<\<F_\alpha)\iso
H^n(f^\times \smalldrlm{\lift 0,\!\alpha{\,},} F_\alpha).
$$
\end{cor}

\emph{Remarks.} 1.~Conditions (i) and (ii) in Theorem~\ref{Thm2.7} are connected via the
 flat, affine, and universally bicontinuous natural map
\mbox{$\Spec(S_{\le1}(F))\to Y\<$}, where $S_{\le 1}(F)$ is the $\SO_Y$-algebra $\SO_Y\oplus F$ with $F^2=0$.\vspace{1pt}

2.~The idea behind the proof of Corollary~\ref{2.7.1} is to use Lazard's theorem that over a commutative ring $A$ any flat module
is a \smash{$\smalldrlm{}\!\!$}\vspace{.8pt} of finite-rank free \dash A.modules. \cite[p.\,163, (6.6.24)]{GD}, to show that (i) and (ii)\vs1 imply
condition~(ii) in Theorem~\ref{Thm2.7}.

\section{Boundedness of $f^\times$ implies finite tor-dimension}

\begin{thm}\label{Thm3.1} Let $f\colon  X\to Y$ be a map.
If\/ $f^\times$ is bounded then $f$~has finite tor\kern.5pt-dimension.

\end{thm}

The proof uses the following two Lemmas.\vs2

An \dash\OX.complex. $E$ is \emph{$a\_\mkern1.5mu$locally projective} $(a\in\Z)$
if there is a $b\ge a$ and an affine open covering 
$\bigl(U_i\set\Spec(A_i)\bigr)_{i\in I}$ of~$X$
such that for each~$i\in I$, the restriction $E|_{U_i}$\
is \dash\mathbf D.isomorphic. to a quasi-coherent direct summand of~a complex~$F$ of free 
\dash\SO_{U_i}.modules., with $F$ vanishing in all degrees outside $[a,b\>]$. 

Every complex with perfect amplitude in $[a,b\>]$ (\S2)  is $a\_\>\>$locally projective.

\begin{lem}\label{Lemma3.3} For any scheme\/ $X,$ there is an integer\/ $s>0$ such that
for all\/~$a\in\Z$ and\/ $a\_\mkern1.5mu$locally projective $E\in\mathbf D(X),$ if\/
$G\in\Dqc(X)$ and\/ $H^jG=0$ for all\/ $j> a-s$ then
$\Hom_{\D(X)}(E,G)=0$.
\end{lem}

\begin{lem}\label{Lemma3.4} Let\/ $f\colon X\to Y$ be a 
perfect map, of tor-dim\/ $d<\infty$. Then there exists an 
integer $t>0$ such that for any $a\_\>$locally projective\/ 
\mbox{$E\in\Dqc(X),$} $\Rf E\in\Dqc(Y)$ is\/ $(a-d-t)\_\>\>$locally projective.
\end{lem}

These Lemmas are proved below. \vspace{1pt}

\medskip

\emph{Proof of~Theorem~\ref{Thm3.1}.} \vs1

Part (i) of Proposition~\ref{Prop2.4} gives an immediate reduction to the case where $Y$ is
affine, say $Y=\Spec(A).$ We need to show in this case that for any open
immersion
$\iota\colon U\hookrightarrow X$
with $U$~affine the \dash\SO_Y.module. $f_*\iota_*\SO_U$ has finite tor-dimension.

Since $U$ is affine, there are natural isomorphisms
$$
f_*\iota_*\SO_U=(f\iota)_*\SO_U\iso\R(f\iota)_*\SO_U\iso
\Rf\R\>\iota_*\SO_U.
$$
So for any $G\in \Dqc(Y)$ there are natural isomorphisms
$$
\Hom_{\mathbf D(Y)}\<(f_*\iota_*\SO_U, G)\cong
\Hom_{\mathbf D(Y)}\<(\Rf\R\iota_*\SO_U, G)\cong
\Hom_{\mathbf D(X)}\<(\R\iota_*\SO_U, f^\times\<G).
$$

Lemma~\ref{Lemma3.4} provides an integer $t\>$ such that
if $\>U$ is\/ \emph{any} quasi-compact open subscheme of $X\<,$ with inclusion 
$\iota\colon U\subset X,$  then\/ 
$\R\>\iota_*\SO_U$ is $(-t)\_\>\>$locally projective.
By Lemma~\ref{Lemma3.3} and the boundedness of $f^\times\<$, it follows,
for $U$ affine, $G$~a quasi-coherent $\SO_Y$-module, and some 
$j\gg0$ not depending on $G$,  that\vspace{-2pt}
\begin{multline*}
\Ext^j(f_*\iota_*\SO_U, G)=\Hom_{\mathbf D(Y)}(f_*\iota_*\SO_U, G[\>j\>])\\
\cong
\Hom_{\mathbf D(X)}(\R\>\iota_*\SO_U, f^\times\<G[\>j\>])= 0.
\end{multline*}
The natural equivalences
$
\mathbf D(A)\overset{\lift.4,\approx,}\lra\mathbf
D(Y_\qc)\overset{\lift.4,\approx,}\lra\Dqc (Y)
$
(where $Y_\qc$ is the category of quasi-coherent $\SO_Y$-modules---see \cite[p.\,30, Cor.\,5.5]{BN}) show then that $f_*\iota_*\SO_U$ has a resolution by the sheafification of a
bounded projective \dash A.complex., and thus has finite tor-dimension, as desired.
\qed

\bigbreak
\stepcounter{thm}
\emph{Proof of~Lemma~\ref{Lemma3.3}.} \vs2

Let us call an  open $U\subset X$ $E$-\emph{good} if $U$ is affine, say $U=\Spec(A)$, and if there is a $b\ge a$ such that the restriction $E|_{U}$\
is \dash\mathbf D.isomorphic. to the sheafification of a projective \dash
A.complex.~E vanishing in all degrees outside $[a,b\>]$.

Clearly, every quasi-compact open subset of $X$ is a finite union of $E$-good open subsets. Hence,  as in the proof of \cite[p.\,13, Prop. 3.3.1]{BB}, it will suffice
to show that Lemma 3.2 holds for $X$ if $X$ itself is $E$-good, or if $X=X_1\cup X_2$ with $X_1$ and $X_2$ quasi-compact open subsets such that Lemma~3.2 holds for $X_1$, $X_2$ and 
$X_1\cap X_2$ (which is also quasi-compact, since $X$ is quasi-separated).

Suppose first that  $X=\Spec(A)$ is $E$-good. Let
E  be as in the definition of $E$-good, and let $G\in\Dqc(X)$ be such that 
$H^jG=0$ for all $j>a-1$.
The natural equivalence of categories
\mbox{$\D(X_\qc)\xrightarrow{\lift.4,\approx,}\Dqc (X)$}  (where
$X_\qc$ is the category of 
quasi-coherent $\OX$-modules) 
allows us to assume  $G$~quasi-coherent, so that 
$G$ is the sheafification of an $A$-complex G\kern1pt; and further, after applying the well-known truncation functor 
(see e.g.,  \cite[\S1.10]{Lp}) we can assume that 
G~vanishes in all degrees $> a-1$.\vadjust{\kern1pt}

The dual versions of \cite[(2.3.4)
and~(2.3.8)(v)]{Lp}, and the  equivalences 
$\D(A)\xrightarrow{\lift.4,\approx,}\D(X_\qc)$,  \mbox{$\D(X_\qc)\xrightarrow{\lift.4,\approx,}\Dqc (X)$},
yield natural  isomorphisms, with $\mathbf K(A)$ the homotopy category of $A$-complexes:
$$
\Hom_{\mathbf K(A)}(\text{E}, \text{G})\cong\Hom_{\mathbf D(A)}(\text{E}, \text{G})\\
\cong \Hom_{\D(X_\qcc)}(E, G)\cong\Hom_{\mathbf D(X)}(E, G).
$$
So since E vanishes in all degrees $<a$ and G vanishes in all degrees $>a-1$, therefore 
$\Hom_{\mathbf D(X)}(E,
G)=0$, proving Lemma 3.2 in this case.\vs1

Suppose next that $X=X_1\cup X_2$ as above. Let $s>0$ be such that Lemma~3.2 holds with this $s$ for all three of $X_1$, $X_2$, and $X_1\cap X_2$.
 Let $G\in\Dqc(X)$ satisfy $H^jG=0$ for all $j>a-(s+1)$.
Let $i\colon X_1\hookrightarrow X$, $j\colon X_2\hookrightarrow X$, and
$k\colon X_1\cap X_2\hookrightarrow X$ be the inclusion maps. One gets
the natural\- triangle
$$
\CD
G@>>> \R i_*i^*\<\<G\oplus\R j_*j^*\<\<G@>>> \R k_*k^*\<\<G@>>> G[1],
\endCD
$$
by applying the usual exact sequence, holding for any flasque
\dash\OX.module.~$F\<$,
$$
0\to F\to  i_* i^*F\oplus j_* j^*F\to
k_*k^*F\to 0
$$
to an injective q-injective resolution%
\footnote{Another term for ``q-injective" is ``K-injective"---see \cite[(2.3.2.3), (2.3.5)]{Lp}.}
of $E^+\<$.
There results an exact sequence \cite[p,\,21, 1.1(b)]{H}, with $\Hom\set\Hom_{\mathbf
D(X)}$,
$$
\Hom(E,\R k_*k^*\<\<G[-1])\to
\Hom(E,G)\to \Hom(E,\R i_*i^*\<\<G)\oplus\Hom(E,\R j_*j^*\<\<G).
$$
Adjointness of $\R k_*$ and $\Ll k^*=k^*$ gives that 
$$
\Hom_{\D(X)}(E,\R k_*k^*\<\<G[-1])\cong \Hom_{\D(X_1\cap
X_2)}(k^*\<\<E,k^*G[-1]);
$$
and Lemma 3.2 makes these groups vanish.
Similarly,
$\Hom(E,\R j_*j^*\<\<G)=0$ and \mbox{$\Hom(E,\R i_*i^*\<\<G)=0$.}
Hence $\Hom(E,G)=0$. 
\qed

\medbreak

\emph{Proof of Lemma~\ref{Lemma3.4}.}

\medskip

The question is local on $Y\<$, so we may assume $Y$ affine, say $Y=\Spec(B)$.

Arguing as in the preceding proof, suppose first that $X$ is $E$-good. We begin with the case $E=\OX$. 
Then for some
\mbox{$t>0$,} $f$~factors as\vs1 
$$
X\xrightarrow{\,\iota\,}Y_t\set\Spec(B[T_1,T_2,\dots,T_t])
\xrightarrow{\pi}\Spec(B),\vs1
$$ 
where $T_1,\dots,T_t$ are independent 
indeterminates, $\iota$ is a closed immersion, and $\pi$ is the 
natural map. 
By \cite[p.\,252, Prop.\,4.4(ii) and p.\,174, Prop.\,2.2.9(b)]{I}, 
the sheaf $\iota_*\OX$ is $\D(Y_n)$-isomorphic to
a bounded quasi-coherent complex~$G$ of direct summands of finite-rank free
$\SO_{Y_n}$-modules, vanishing in all degrees \mbox{$<-\>d-\<t$.} Hence
$\Rf\OX\cong \pi_*\iota_*\OX \cong \pi_*G$ is 
$(-\>d-\<t)\_\>\>$locally projective.

Since $\Rf$ commutes with direct sums in $\Dqc$ (because $\Rf$ has a
right adjoint, or more directly, by \cite[3.9.3.3]{Lp}), it follows that for any free
$\OX$-module~$E$, $\Rf E$ is 
$(-\>d-\<t)\<\_\>\>$locally projective. Finally, to show that for any
$a\_\>\>$locally projective $E$, $\Rf E$ is 
$(a-d-t)\<\_\>\>$locally projective, one reduces easily to where
$E$ is a bounded free complex, and then argues by induction on the
number of degrees in which $E$ is nonvanishing, using the following
observation:

\smallskip

\noindent $(*)\colon$ 
\emph{In a\/ $\D(X)$-triangle\/ $N[-1]\xrightarrow{\delta}L\to M\xrightarrow{\rho}N,$
if\/ $N$ is $a\_\>\>$locally projective \phantom{$\!(**)\colon$}then $M$ is 
$a\_\>\>$locally projective  iff so is} $L$.
 
 \smallskip
(To see this, one may suppose that $X$ is affine, say $X=\Spec(A)$.
If~$N$ and~$L$ are  $a\_\>\>$locally projective, then one may assume they are
 sheafifications of  bounded projective
\dash A.complexes. vanishing in all degrees {$<\<\<a$}, so that  by
the dual versions of \cite[(2.3.4)
and~(2.3.8)(v)]{Lp},  $\delta$ comes from
a \dash\mathbf K(X).morphism. $\delta_0\colon N[-1]\to L\>$; and $M$ is
isomorphic to the mapping cone of $\delta_0$, an $a\_\>\>$projective
complex. Similarly, if $N$ and $M$ are   sheafifications of  
bounded projective
\dash A.complexes. vanishing in all degrees {$<\<a$}, then $\rho$ comes from 
a \dash\mathbf K(X).morphism. $\rho^{}_0\colon M\to N$, and since $L[1]$  is isomorphic to the mapping cone of~$\rho^{}_0$, therefore $L$ is $a\_\>\>$locally projective.)\vs3\looseness=2

Suppose next that $X=X_1\cup X_2$ with $X_1$ and $X_2$ quasi-compact open subsets for which there exists a $t>0$  such that Lemma~3.3 holds with this $t$ for all three of $X_1$, $X_2$, and $X_1\cap X_2$.
As in the proof of Lemma~\ref{Lemma3.3}, 
there is a $\D(Y)$-triangle
$$
\minCDarrowwidth=.17in
\CD
\Rf E\<@>>> \<\R(fi)_*(E|_{X_1})\<\oplus\<\R(fj)_*(E|_{X_2})\<@>>>\<\R(fk)_*(E|_{X_1\cap X_2})
\<@>>>\< \Rf E[1]
\endCD
$$
in which
the two vertices 
other than $\Rf E$ are $(a-d-t)\_\>\>$projective, whence, by ~$(*)$, so is $\Rf E$. 
\vspace{3pt}

As before, this completes the proof of Lemma~\ref{Lemma3.4}, and so of Theorem~\ref{Thm3.1}.

\section{Approximation by perfect complexes.}\label{S0}

Terminology remains as in \S2.\vs1

The main results in this section are the following two theorems. 

\thmm{perfapprox} For any\/ scheme $X\<,$ 
there exists a positive integer $B=B(X)$
such that for any\/ $E\in\Dqc(X)$ and integer\/~$m,$  
if\/ $E$ is\/ $(m-B)$-pseudo-coherent then
there exists in\/ $\Dqc(X)$ an\/ \mbox{$m$-isomorphism\/} $P\to E$ with $P$ perfect.
\ethm

 \thmm{perfgen}
Let $X$ be a scheme. Then\/ $\Dqc(X)$ has a perfect generator, i.e., there is a perfect\/ $\OX$-complex\/ $S$ such that for each\/ 
$E\ne 0$ in\/ $\Dqc(X)$ there is an\/ $n\in\mathbb Z$ and a nonzero 
$\Dqc(X)$-morphism\/ $S[n]\to E$.

Moreover, for each such\/ $S$ there is  an integer\/ $A=A(S)$ such that
for all\/~$E\in\Dqc(X)$ and $j\in\mathbb Z$ with\/ $H^j\<(E)\neq0,$ 
\[
\Hom\big(S[n], E\big)\neq0\quad\textup{for some }\,n\le A- j.\vs1
\]
\ethm

Theorem~\ref{perfapprox} may be compared to \cite[p.\,173, 2.2.8(b)]{I}). The first
statement in Theorem~\ref{perfgen} comes from \cite[p.\,9, Thm.\,3.1.1]{BB}.

Proofs are given in section 5 below.

\stepcounter{thm}

\begin{cor}\label{qp bounded}
If a map\/ 
$\>f$ is either perfect or quasi-perfect, then the functor\-~$f^\times\!$ is bounded.
\end{cor}

\prf
As mentioned in the Introduction,
$f^\times$~commutes with translation of complexes, and
$f^\times$ is bounded below. So to show that $f^\times$ is bounded, 
it is enough to find a~ $j^{}_0$
such that for every $m\in\mathbb Z$ and $F\in\Dqc(Y)$ with $H^i(F)=0$ for all $i>m$, 
it holds that $H^j\<f^\times\<\< F =0$ for all $j\ge m+j^{}_0\>$.

Suppose $H^{j}(f^\times F)\ne 0$. With $S$ and $A$ as in Theorem~\ref{perfgen},
there exists \mbox{$k\le A$} and a nonzero $\mathbf D(X)$-morphism 
$S\to f^\times\<\< F[\>j-k]$, the latter corresponding under adjunction to  a \emph{nonzero} morphism 
$\lambda\colon\Rf S\to F[\>j-k]$.

For some $a$, $\Rf S$ is $a\_\>\>$locally projective---when $f$ is perfect, that results from Lemma~\ref{Lemma3.4}, and when $f$ is quasi-perfect, it's because $\Rf S$ is perfect.
It follows from Lemma~\ref{Lemma3.3} that
there is an integer $s=s(Y)$ such that $\lambda$~cannot exist if~$\>j\ge m+A-a+s$.
With $j^{}_0\set A-a+s$, we must have then that 
$H^{j}(f^\times F)=0$ for all $\>j\ge m+j^{}_0\>$; 
and  so $f^\times$ is  indeed bounded.
\eprf

\begin{cor}\label{perf->pc}
For a map\/ $f\colon X\to Y,$ the following are equivalent.

{\rm (i)} $f$ is quasi-proper. 

{\rm(ii)} For any
perfect\/ $\OX$-complex\/~$P,$ $\Rf P$ is pseudo-coherent.

{\rm(iii)} If\/ $S$ is as in Theorem~\ref{perfgen}, then\/ $\Rf S$ is pseudo-coherent.
\end{cor}

\prf
(i)${}\Rightarrow{}$(ii)${}\Rightarrow{}$(iii). The first implication is clear (since perfect complexes are 
pseudo\kern.5pt-coherent); and the second is trivial. \vs1

(iii)${}\Rightarrow{}$(ii).  Let $R$ be the smallest triangulated\vs{-1.5} subcategory of 
$\>\Dqc(X)$ containing $S$, and let $\widehat R$ be the full subcategory\vs{-1.5}   of $\Dqc(X)$
whose objects are  all the direct summands of objects
in $R$. 
The subcategory $\widehat R\subset\Dqc(X)$ is triangulated, and closed under  formation of direct summands 
\cite[p.\,99, 2.1.39]{N2}.\looseness=-1

\enlargethispage*{3pt}
 The full subcategory $R^c$ of $R$ whose objects are the compact ones in~$R$ is triangulated, whence\vs{.3} every object in $R$---and in $\widehat R$---is compact.
Consequently, \vs{-1} \cite[p.\,222, Lemma 3.2]{N1} shows that the smallest full subcategory of~ $\Dqc(X)$ which contains~$\widehat R$ and is closed with respect to coproducts is $\Dqc(X)$ itself. 

\noindent Hence, by
\cite[p.\,214, Thm.\,2.1.3]{N1}, every perfect complex lies in $\widehat R$. 
(Alternatively, see \cite[p.\,285, Prop. 8.4.1 and p.\,140, Lemma 4.4.5]{N2}.)\vs1

Since the pseudo\kern.5pt-coherent complexes in $\Dqc(Y)$ are the objects of a triangulated subcategory closed under formation of direct summands 
\cite[p.\,99,~b) and p.\,105, 2.12]{I}, therefore the complexes
$Q\in\Dqc(X)$ such that $\Rf Q$ is pseudo\kern.5pt-coherent are the objects of a triangulated subcategory closed under formation of direct summands. Thus if $S$ is such a $Q$ then every complex
in~$\widehat R$---and so every perfect complex---is such a $Q$. \vs1

(ii)${}\Rightarrow{}$(i).
Let $E$ be a pseudo\kern.5pt-coherent $\OX$-complex, let $m\in\mathbb Z$, and
let 
$$
\CD
P@>\alpha>>E@>>>Q@>>>P[1]
\endCD
$$
be a triangle with $\alpha$ an $m$-isomorphism as in Theorem~\ref{perfapprox}. Thus \mbox{$H^k(Q)=0$} for all $k\ge m$. As $\Rf$ is 
bounded above \cite[(3.9.2)]{Lp}, there is an integer $t$ depending only on~$f$
such that $H^k(\Rf Q)=0$ for all $k\ge m+t$,  that is, $\Rf\alpha$ is an
$(m+t)$-quasi-isomorphism. So if $\>\Rf P$ is
pseudo\kern.5pt-coherent then $\Rf E$ is $(m+t)$-pseudo\kern.5pt-coherent; and
since $m$ is arbitrary, therefore $\Rf E$ is pseudo\kern.5pt-coherent.
\eprf

From \ref{perf->pc}(ii) we get:

\begin{cor}\label{qpqp} Every quasi-perfect map is quasi-proper.
\end{cor}

\smallskip
Next, we deduce  ``stability" of quasi-properness.

\begin{prop}\label{stableqp}
Let
$$
\begin{CD}
X'@>v>>X  \\
@V g VV  @VV f V  \\
Y' @>>\lift1.2,u,> Y
\end{CD}
$$
be a tor-independent square. If  $f$ is
quasi\kern.5pt-proper then so is $g$.
\end{prop}

\emph{Proof.}
Since pseudo\kern.5pt-coherence is a local property, it suffices to prove the Proposition when $Y'$ is affine and $u(Y')$ is contained in an affine subset
of~$Y$. So we can assume that $u=u'u''$ where $u'$ is an open immersion and $u''$ is an affine map. It follows that it suffices to prove  the Proposition 
(a) when $u$---hence $v$---is an open immersion and (b) when $u$---hence $v$---is an affine map
(see \cite[p.\,358, (9.1.16)(iii), (9.1.17)]{GD}). 

In either of these two cases, it holds that

\noindent $(*)$ \emph{ if\/ $S$ is as in 
Theorem~\ref{perfgen} then\/ $\Ll v^*\<\<S$ is a generator of} $\Dqc(X')$.\vs1

Indeed, in case $v$ is an open immersion and $0\ne E\in\Dqc(X')$ then
\mbox{$0\ne\R v_*E\in\Dqc(X)$} (since $E\cong v^*\R v_*E$); and the same holds in case $v$ is affine, by  \cite[(3.10.2.2]{Lp}. So in either case,  for some~$n$,
$$
0\ne \Hom_{\D_\qcc(X)}\<(S[n], \R v_*E)\cong \Hom_{\D_\qcc(X')}\<(\Ll v^*S[n], E),
$$
proving $(*)$.

It is easy to see that the complex 
$\Ll v^*\<\<S$ is perfect.
So by Corollary\-~\ref{perf->pc}, to prove the Proposition for $u$ as in $(*)$ it suffices to show that $\R g_*\Ll v^*\<\<S$ is pseudo\kern.5pt-coherent. 
But by \cite[(3.10.3)]{Lp}, $\R g_*\Ll v^*\<\<S\cong \Ll u^*\Rf S$; and since $\Rf S$ is pseudo\kern.5pt-coherent, therefore,
by \cite[p.\,111, 2.16.1]{I}, so is $\Ll u^*\Rf S$.\vs1

\section{Proofs of  Theorems 4.1 and 4.2}\label{proofs}
Heavy use will be made of the following technical notion.

\begin{dfn}\label{D1.0.1}
Let $\ct$ be a triangulated category, and  let $\cs\subset\ct$ be a class of objects.
Let $m\leq n$ be integers. The full subcategory $\cs[m,n]\subset\ct$ is the
smallest among (\,= intersection of) all full subcategories $\mathcal S\subset\ct$ such that:\vs{-1}
\begin{enumerate}

\item[(i)] 
$0$ is contained in $\mathcal S$.

\item[(ii)] If $E\in \cs$, then $E[-\ell\>]\in \mathcal S$ for all integers $\ell$ in the interval $[m,n]$.

\item[(iii)]
For any   $\ct$-triangle\vs{-2}
\[
\CD
E @>>>F @>>> G @>>>  E[1]\>,\vs{-2.5}
\endCD 
\]
if $E$ and $G$ are in $\mathcal S$ then so is $F$.
\end{enumerate}
\edfn 

\rmk{R.1.0.2} One checks that 
$\cs[m,n]=\big(\bigcup_{\ell=m}^n \cs[-\ell\>]\>\big)\<[0,0]$.\vs1
\ermk

\rmk{R1.0.2}Defn.\,\ref{D1.0.1} expands\vs1
to allow  $m=-\infty$ or $n=\infty$. For example,
$\cs[m,\infty)\set \bigcup_{n=m}^\infty \cs[m,n]$. Furthermore,\vs{1.2}
$\cs(\infty,\infty)\set \bigcup_{\lift1.3,m\le n,} \cs[m,n]$, being closed under translation
(see ~\ref{R1.0.3}(i)), is the smallest triangulated\vs{.3} subcategory of\/ $\ct$ 
containing\/ $\cs$ \cite[p.\,60, Defn.\,1.5.1]{N2}.
\ermk

\rmk{R1.0.3}
The following are easy observations.
\begin{enumerate}
\item[(i)]
\emph{If\/ $E\in \cs[m,n]$ and\/ $j\in\mathbb Z\>$ then\/ $E[-j\>]\in \cs[m+j,\>n+j\>]$.}

\noindent Indeed,  (i), (ii) and (iii) in \ref{D1.0.1}  hold
for the full subcategory $\mathcal S\subset \ct$ whose~objects are those $E\in\cs[m,n]$ such that 
$E[-j\>]\in \cs[m+j,\>n+j\>]$.\looseness=-1

\noindent One deduces that, with  $\cs[m,n]_{\rm o}$ the class of objects in $\cs[m,n]$,
$$
\big(\cs[m,n]_{\rm o}\big)[m'\<,\>n'\>]=\cs[m+m'\<,\>\>n+n'\>].
$$

\item[(ii)]
\emph{If every object of\/ $\cs$ is compact, then so is every object
of} $\cs[m,n]$.\vs1

\noindent  Indeed,  (i), (ii) and (iii) in \ref{D1.0.1}  hold
for the full subcategory $\mathcal S\subset \ct$ whose~objects are those $E\in\cs[m,n]$ which are compact.\vs1

\item[(iii)]
\emph{Let\/ $\ca$ be an abelian category, and\/ $H:\ct\la\ca$  a cohomological functor, 
see\/}  \cite[p.\,32, 1.1.9]{N2}. \emph{If for every object $F\in \cS$
we have $H(F[-i\>])=0$ for all\/ $i$ in some interval\/ $[a,b\>],$ then\/
for all} $\>E\in \cs[m,n]$, 
\begin{equation*}\label{van}
H(E[-j\>])=0\> \textup{ for  all } j\in[a-m,\>b-n].\tag{\ref{R1.0.3}.1}
\end{equation*}

\noindent  Indeed,  (i), (ii) and (iii) in \ref{D1.0.1}  hold
for the full subcategory $\mathcal S\subset \ct$ whose~objects are those $E\in\cs[m,n]$ which satisfy \eqref{van}.\vs1

\item[(iv)] \emph{Let\/ $\phi\colon\ct\to\ct'$ be a triangle-preserving additive functor} \cite[\S1.5]{Lp}.
\emph{Then} 
$$
\phi\big(\cs[m,n]\big)\subset \big\{\phi\cs\big\}[m,n].
$$
Indeed,  (i), (ii) and (iii) in \ref{D1.0.1}  hold
for the full subcategory $\mathcal S\subset \ct$ whose~objects are those $E\in\cs[m,n]$ such that $\phi E\in \big\{\phi\cs\big\}[m,n]$.
\item[(v)]
\emph{Let\/ $A\xrightarrow{\alpha}B\xrightarrow{\beta} C$ be\/ morphisms inside\/  $\ct$-triangles}
\[
\CD
E @>>> A @>\alpha>> B @>>> E[1]@., \\
F @>>> A @>\beta\alpha>> C @>>> F[1]
@., \\
G @>>> B @>\beta>> C @>>> G[1] @..\\[3pt]
\endCD
\]
\emph{If\/ $E$ and\/ $G$ are in\/ $\cs[m,n]$ then so is} $F$.\vs{2}

\noindent Indeed, the octahedral axiom \cite[p.\,60, 1.4.7\kern.5pt]{N2}
produces a   triangle\vs{-1} 
\[
\CD
E@>>> F@>>>G@>>> E[1].
\endCD
\]

\end{enumerate}
\ermk

\exm{E1.0.5} 
\textup{Remark~\ref{R1.0.3}(iii) will be used thus. 
Let $G$ be an object of $\ct$, and $H$ the cohomological functor $H(-)\set\Hom(-,G)$,
see \cite[p.\,33, 1.1.11]{N2}.
Then for $a=m$ and $b=n$ the assertion becomes:}\vs1

If for every object $F\in \cs$
we have $\Hom(F[-i\>],G)=0\>$ for all\/ $i\in[m,n],$ then\/ $\Hom(E,G)=0$
for all\/ $E\in \cs[m,n]$.
\eexm

 \smallskip
 
 A key role in the proofs will be played by \emph{Koszul complexes.}

\begin{ex}\label{e1.0.5}
Let $R$ be a commutative ring,  $(f_1,f_2, \dots,f_r)$
a sequence in~$R$, and $(n_1,n_2, \dots,n_r)$ a sequence of positive integers. 
The associated\vs{1.2} Koszul complex (see, e.g.,  \cite[III, (1.1.1)]{EGA}) is\vs{-1.5}
$$
K_\bullet(f_1^{n_1}\<,\dots,\>f_r^{n_r})\set \otimes_{i=1}^r K_\bullet(f_i^{n_i}),\\[3pt]
$$
where  $K_\bullet(f_i^{n_i})$ is $R\xrightarrow{\lift1.3,f_i^n,}R$ 
in degrees $\>-1$ and 0,
and $(0)$ elsewhere. 
 
For $r=0$, set $K_\bullet(\phi)\set R$. For all $r\ge 0$, $K_\bullet(f_1^{n_1},\dots,f_r^{n_r})$ is a complex with perfect amplitude in $[-r,0]$, and
 homology killed by each $f_i^{n_i}\<$.\vs1

For any complex $E$, and $f\in R$, $K_\bullet(f)\otimes E$ is the mapping cone of
 the endomorphism
``multiplication by $f\>$'' of $E$. Thus for $1\le i<r$, 
$K_\bullet(f_i^{n_i},\dots,f_r^{n_r})$ is the mapping cone 
of the endomorphism ``multiplication by~$f_i^{n_i}\>$'' of the complex 
$K_\bullet(f_{i+1}^{n_{i+1}}\<,\>f_{i+2}^{n_{t+2}},\dots,\>f_r^{n_r})$.
It follows that$$
\>K_\bullet(f_1^{n_1}\<,\>f_2^{n_2},\dots,\>f_r^{n_r})\in
\left\{K_\bullet(f_1,f_2,\dots,\>f_r)\right\}[0,0].
$$
This is shown by a straightforward induction, based on application of~\ref{R1.0.3}(v)
to the following  three natural triangles (where $\widehat{\phantom{m}}$ signifies ``omit,"):\vs{-2}
\begin{multline*}
K_\bullet(f_1^{n_1}\<,..\>,\>f_i^{n_i}\<\<,\>f_{i+1},..\>,\>f_r) \xrightarrow{\quad\: }
\\[-5pt]
  K_\bullet(f_1^{n_1 }\<,..\>,\>\widehat{f_i^{n_i}},f_{i+1},..\>,\> f_r)[1]
\xrightarrow{ f_i^{n_i}}  K_\bullet(f_1^{n_1}\<,..\>,\>\widehat{f_i^{n_i}},\>f_{i+1},..\>,\>f_r)[1]
  \xrightarrow{\ +\ }\mkern-15mu
\end{multline*}
\vs{-12}
\begin{multline*}
K_\bullet(f_1^{n_1}\<,..\>,\>f_i^{n_i^{}\<+1}\<,\>f_{i+1},..\>,\>f_r) \xrightarrow{\quad\: }\\[-6pt]
\< K_\bullet(f_1^{n_1}\<,..\>,\>
 \>\widehat{\<\strut  f_i^{{\vbox to 3pt{\vskip-.5 ex\hbox{$\scriptstyle n_i^{}\<+1$}\vss}}}\,\,}\!\!,
    \>f_{i+1}\<,..\>,\>f_r)[1]
\xrightarrow{\! f_i^{n_i^{}\<+1}\!\!\!}  K_\bullet(f_1^{n_1}\<,..\>,\>
\>\widehat{\<\strut  f_i^{{\vbox to 3pt{\vskip-.5 ex\hbox{$\scriptstyle n_i^{}\<+1$}\vss}}}\,\,}\!\!,
\>f_{i+1},..\>,\>f_r)[1]
    \xrightarrow{\ +\ }\mkern-15mu
\end{multline*}
\vs{-12}
\begin{multline*}
K_\bullet(f_1^{n_1}\<,..\>,\>f_i,\>f_{i+1},..\>,\>f_r) \xrightarrow{\quad\: }\\[-5pt]
K_\bullet(f_1^{n_1}\<,..\>,\>\widehat{f_i},\>f_{i+1}\<,..\>,\>f_r)[1]
\xrightarrow{\  f_i\ }  K_\bullet(f_1^{n_1}\<,..\>,\>\widehat{f_i},\>f_{i+1},..\>,\>f_r)[1]
   \xrightarrow{\ +\ }\mkern-15mu
\end{multline*}
\end{ex}

\medskip

The proofs of Theorems~\ref{perfapprox} and~\ref{perfgen} will 
involve induction on the number of
affine open subschemes needed to cover $X\<$. One needs to begin
with some results on affine schemes. 

In the situation of Example~\ref{e1.0.5},
denote the sequence\vs{.6} $(f_1^n\<\<,\dots,\>f_r^n)\ (n>0)$ by $\>\mathbf f^n\<\<,\>$ omitting the superscript ``$n$" when $n=1$. Let
$C_\bullet(\mathbf f^n)$ be the cokernel of that map of complexes
$R[-1]\to K_\bullet(\mathbf f^n)[-1]$ which is the identity map of~$R\>$ in degree~1.
The complex $C_\bullet(\mathbf f^n)$ has
perfect amplitude in~$[1-r,\>0]$; 
and there is a natural homotopy triangle
\begin{equation*}\label{htri}
\CD
C_\bullet(\mathbf f^n)@>>> R@>>> K_\bullet(\mathbf f^n)@>>> C_\bullet(\mathbf f^n)[1]\>.
\endCD\tag{\ref{e1.0.5}.1}
\end{equation*}

There is a map of complexes $K_\bullet(f^{n+m})\to K_\bullet(f^n)\ (f\in R;\;m,n>0)$ depicted by
$$
\CD
R @= R\\
@Af^{n+m}AA @AAf^nA \\
R@>>f^m> R
\endCD
$$
Tensoring such maps gives a map  
$K_\bullet(\mathbf f^{n+m})\to K_\bullet(\mathbf f^n)$,
and hence a map
$C_\bullet(\mathbf f^{n+m})\to C_\bullet(\mathbf f^n)$.
For any \dash R.complex.~$E$,\vadjust{\kern.6pt}
we have then the \emph{\v Cech complex}
$$
\check C^\bullet(\mathbf f, E)\set
\smalldrlm{\lift-.1,\<n,}\> \Hom^\bullet_R(C_\bullet(\mathbf f^n),\>E).\\[-1pt]
$$

Let  $U$ be the complement of the closed subscheme
$\Spec (R/\mathbf fR)\subset\Spec(R)$, 
with inclusion $\iota\colon U\hookrightarrow X$.
From \cite[III, \S1.3]{EGA}
 it follows readily that, with~$E^\sim$
the quasi-coherent complex corresponding to $E,$  
 $\R\>\iota_*\iota^{\<*}\<\<E^\sim$ \emph{is naturally 
\dash\mathbf
D.isomorphic. to the sheafified \mbox{\v Cech complex\/} 
$\>\check{\mathcal  C}^\bullet(\mathbf f,
E)\set \check C^\bullet(\mathbf f,E)^\sim\<.$}

In particular, if the homology of $E$ is $\mathbf fR$-torsion (i.e., for
all~$i\in \mathbb Z$, each element
of $H^i(E)$ is annihilated by a power\vs{.6} of $\mathbf fR\>$)---or
equivalently, if  $\iota^{\<*}\<\<E^\sim$ is exact---then $\check C^\bullet(\mathbf f,E)$ is exact;
and since the complex $C_\bullet(\mathbf f^n)$ is bounded and  projective, therefore
$$
H^0\Hom^\bullet_R(C_\bullet(\mathbf f^n),\>E)\cong
H^0\R \Hom^\bullet_R(C_\bullet(\mathbf f^n),\>E)\cong
\Hom_{\mathbf D(R)}(C_\bullet(\mathbf f^n),\>E),
$$
whence 
$$
\smalldrlm{\lift-.1,\<n,}\> \Hom_{\mathbf D(R)}(C_\bullet(\mathbf f^n),\>E)
\cong H^0\check C^\bullet(\mathbf f, E)=0.
$$
Consequently, the commutative diagrams of the following form,  with exact rows coming
from ~\eqref{htri}, and  columns  from the maps described above:
$$
\minCDarrowwidth=.25in
\CD
\Hom_{\mathbf D(R)}(K_\bullet(\mathbf f),E)@>>>\Hom_{\mathbf D(R)}(R,E)
@>>>\Hom_{\mathbf D(R)}(C_\bullet(\mathbf f),E)\\
@VVV @| @VVV\\
\Hom_{\mathbf D(R)}(K_\bullet(\mathbf f^n),E)@>>>\Hom_{\mathbf D(R)}(R,E)
@>>>\Hom_{\mathbf D(R)}(C_\bullet(\mathbf f^n),E)
\endCD
$$
show that \emph{for any\/ $\lambda\in\Hom_{\mathbf D(R)}(R,E)$ 
there is an\/ $n>0$ such that\/ $\lambda$~factors through a\/ 
$\mathbf D(R)$-morphism} $K_\bullet(\mathbf f^n)\to E$. \vs1 

If $P$ is a bounded complex of finitely generated projective $R$-modules, then
the homology of $\Hom^\bullet_R(P,E)$ is still $\mathbf fR$-torsion 
(as one sees, 
e.g., by induction on the number of nonvanishing components of $P$); 
and replacing
 $E$ in what precedes by $\Hom^\bullet_R(P,E)$, one obtains, 
via $\Hom$-$\otimes$ adjunction, that 
\emph{for any\/ $\lambda\in\Hom_{\mathbf D(R)}(P,E)$ 
there is an\/ $n>0$ such that\/ $\lambda$~factors\vs{.6} through a\/ 
$\mathbf D(R)$-morphism} $K_\bullet(\mathbf f^n)\otimes P\to E$.

\begin{lem}\label{1.2.0}
 Let\/ $E$ be an\/ $R$-complex such that 
$H^j(E)$ is $\mathbf fR$-torsion for all\/~$j\ge-r,$  $P$ an 
$R$-complex with perfect amplitude in\/ $[0,b\>]$ for some\/ $b\ge0,$ 
and\/ $\lambda\in\Hom_{\mathbf D(R)}(P,E)$. Then there is an integer\/
\/$n>0$ and a\/ homomorphism of\/ $R$-complexes\/ 
\mbox{$\>\lambda_n\colon K_\bullet(\mathbf f^n)\otimes P\to E$} such that  
for all\/ $j\ge-r,$ the homology map
$
H^j(\lambda)\colon H^j(P)\to H^j(E) 
$ 
factors as \vs{-2}
$$
H^j(P)=H^j(R\otimes P)\xrightarrow{\<\textup{natural}\<\<}H^j\big(K_\bullet(\mathbf f^n)\otimes P\big)\
\xrightarrow{H^j(\lambda_n)}H^j(E).
$$
\end{lem}

\begin{proof} We may assume that $P$ is a complex of finitely generated projective $R$-modules, vanishing in all degrees outside $[0,b\>]$, see \cite[p.\,175, b)]{I}. Let\/ $\tg{-r}\> E\>$~be the usual truncation, and $\pi\colon E\to\tg{-r}\> E$ the natural map, which induces homology\vs{-.5} isomorphisms in all
degrees $\ge -r$ (see, e.g., \cite[\S1.10]{Lp}).
By the preceding remarks, $\pi\lambda$ factors in $\D(R)$ as
 $$
P=R\otimes P\xrightarrow{\<\textup{natural}\<\<}
K_\bullet(\mathbf f^n)\otimes P
\xrightarrow{\!\!\!\bar{\ \; \mathstrut\lambda_n}}\tg{-r}\>E.
$$
Since $K_\bullet(\mathbf f^n)\otimes P$ is bounded and projective, we may assume that $\bar\lambda_n$ is a map of $R$-complexes.
Then the $R$-homomorphism 
$$
(\bar\lambda_n)^{-r}\colon 
\big( K_\bullet(\mathbf f^n)\otimes P\big)^{-r}\<=P^0\to (\tg{-r}\>E)^{-r}= \coker(E^{-r-1}\to E^{-r})
$$
lifts to a map $P^0\to E^{-r}\<$, \kern.5pt giving a map $\lambda_n$ with the desired properties.
\end{proof}

\begin{cor}\label{L1.2.1}
Set $I\set\mathbf fR=(f_1,f_2,\dots,f_r)R$. Let\/ $m\in\mathbb Z$ and let\/ $E$ be an\/ $R$-complex such that\/ $H^i(E)$ is\/ $I$-torsion 
for all\/ $i\ge m-r.$ \vs1

{\rm (i)} If\/ $E$ is\/ $m$-pseudocoherent, and\/ 
$p\ge m$ is such that \/$H^i(E)=0\>$ for all\/~$i> p,$
then there exists in the homotopy category of\/ $R$-complexes 
an\/ \mbox{$m$-quasi-isomorphism} $P\to E$
with $P\in \{K_\bullet(\mathbf f)\}[m,\>p\>]$. \vs1

{\rm (ii)} For any\/ $i\ge m$ with\/ $H^i(E)\ne 0$, there is a nonzero map\/
$
K_\bullet(\mathbf f)[-i\>]\to E.
$ 

\end{cor}

\prf
(i) By \cite[p.\,103, 2.10(b)]{I}, $H^p(E)$ is a finitely generated
$R$-module. So there is an  $\ell>0$ and a surjective homomorphism
$R^\ell\twoheadrightarrow H^p(E)$,  which lifts
to $R^\ell\to \ker(E^p\to E^{p+1})$, and thus there is a homomorphism 
$R^\ell\to E[p]$, or equivalently, $\lambda\colon  R^\ell[m-p\>]\to E[m]$, giving rise, by 
Lemma~\ref{1.2.0}, to an $R$-homomorphism 
$$
\lambda_n[-m]\colon P^{}_1\set \big(K_\bullet(\mathbf f^n)\otimes R^\ell[-p\>]\big) \to E
$$ 
such that $H^{p}\big(\lambda_n[-m]\big)$ is surjective.\vs{.6}
By Example~\ref{e1.0.5} and \mbox{Remark~\ref{R1.0.3}(i),} we have $K_\bullet(\mathbf f^n)[-p\>]\in \{K_\bullet(\mathbf f)\}[p,\>p\>]$.
So we get a homotopy triangle 
\[
\CD
P^{}_1 @>>> E @>\alpha>> Q^{}_1 @>>> P^{}_1[1]
\endCD
\]
with $P^{}_1\in \{K_\bullet(\mathbf f)\}[p,\>p\>]$
and $H^i(Q^{}_1)=0$ for \mbox{all $i\geq p$,}\vs{.6} giving (i) when $p=m$.\looseness=-1  

In any case, $Q^{}_1$ is $m$-pseudocoherent \cite[p.100, 2.6]{I}; and since all the homology of $P^{}_1$ is $I$-torsion, 
the exact homology sequence of the preceding triangle shows that $H^i(Q^{}_1)$ 
is\/ $I$-torsion for all\/ $i\ge m-r$. 
If $m<p$ then, using induction
on $p-m$, one may assume that there is a homotopy triangle\looseness=-1
\[
\CD
P^{}_2 @>>> Q^{}_1 @>\beta>> Q @>>> P^{}_2[1]
\endCD
\]
with $P^{}_2\in \{K_\bullet(\mathbf f)\}[m,\>p-1]\subset  \{K_\bullet(\mathbf f)\}[m,\>p\>]$ 
and\vs{.6}  $H^i(Q)=0$ for all $i\geq m$. 

\noindent There exists then a homotopy triangle 
\[
\CD
P @>>>E @>\beta\alpha>> Q @>>> P[1]
\endCD
\]
which, by Remark~\ref{R1.0.3}(v), is as desired.\vs{1.5}

(ii) There is, by assumption, a nonzero map $R\to H^i(E)$, which lifts
to a map $R\to \ker(E^i\to E^{i+1})$; and so there is a nonzero map
$\lambda\colon R\to E[i\>]$ with $H^0(\lambda)\ne 0$.
If~$j\ge -r$ then $j+i\ge i-r\ge m-r$, so $H^j(E[i\>])=H^{j+i}(E)$ is $I$-torsion, whence by Lemma~\ref{1.2.0}, there is for some  $n>0$ a nonzero map $K_\bullet(\mathbf f^n)\to E[i\>]$.
By Example~\ref{e1.0.5}, $K_\bullet(\mathbf f^n)\in K_\bullet(\mathbf f)[0,0]$; so by 
Example~\ref{E1.0.5}, there is a nonzero map $K_\bullet(\mathbf f)\to E[i\>]$, proving (ii).
\eprf

For dealing with the nonaffine situation, we need to set up some notation.
\ntn{N1.3}
A  scheme 
$X$ can be covered
by finitely many open affine subsets, say $X=\bigcup_{k=1}^tU_k$, with 
\mbox{$U_k =\Spec(R_k)$.} 
For $1\leq k\leq t$, set 
\begin{enumerate}
\item[(i)]
$V_k\set\bigcup_{i=k}^tU_i$.\vs1
\item[(ii)]
$Y_k\set X-V_{k+1}$ ($\>\set X$ when $k=t)$.\vs1
\end{enumerate}

So we have a filtration by closed subschemes 
$
Y_1\subset Y_2\subset\dots\subset Y_t=X.
$

Both $U_k$ and $V_{k+1}$ are quasi-compact open subsets of the (quasi-separated) scheme~$X$, whence so is $U_k\cap V_{k+1}$.
So there is a sequence
$$
\mathbf f_k=\{f_{k1},f_{k2},\cdots,f_{kr^{}_{\<\<k}}\} 
$$
in~$R_k$ such that\vs{-3}
\[
U_k\cap V_{k+1}\:=\:
\bigcup_{i=1}^{r^{}_{\<\<k}}\hbox{Spec}(R_k[1/f_{ki}\>])\ .
\]

\pagebreak[3]

Set
\begin{enumerate}
\item[(iii)]
$I_k\set\mathbf f_kR_k$, 
(so that $U_k\cap V_{k+1}$ is the complement  of 
the closed subscheme $\Spec(R_k/I_k)\subset U_k$). 

\item(iv)
$C_k\set \bigl(K_\bullet(\mathbf f_k)\oplus K_\bullet(\mathbf f_k)[1]\bigr)^\sim=
\bigl(K_\bullet(0,f_{k1},f_{k2},\cdots,f_{kr^{}_{\<\<k}})\bigr)^\sim$\vs{1.5} \newline
with $K_\bullet(*)$  the Koszul complex over $R_k$ associated to 
the sequence~$(*)$, and 
 $(-)^\sim$  the sheafification functor from $R_k$-modules to quasi-coherent
$\SO_{U_k}$-modules---so that $C_k$ is a perfect $\SO_{U_k}$-complex.\vs1 

(The reason for introducing this $\oplus$ will emerge shortly.)
 \end{enumerate}

We have the cartesian diagram of (open) inclusion maps
$$
\CD
U_k\cap V_{k+1} @>\nu>> V_{k+1} \\
@V\lambda VV @VV\xi V \\
 U_k@>>\lift1.1,\mu,> V_k=\!\quad \makebox[0pt]{\quad $U_k\cup V_{k+1}$\hss}
 \endCD
 $$
 
The restriction $\lambda^*C_k$ is homotopically
trivial, whence, in $\mathbf D(V_{k+1})$,
$$
\xi^*\R\mu_*C_k\cong\R\nu_*\lambda^*C_k=0.
$$
Thus, the restrictions of $\R\mu_*C_k$ to both $V_{k+1}$ and $U_k$ are
perfect, and so $\R\mu_*C_k$~is itself perfect.\vs1

For any $\SO_{V_k}$-complex  $G$, the obvious triangle
\[
\CD
G@>0>>G@>>> G\oplus G[1]@>>> G[1]
\endCD
\]
shows that  the complex $G\oplus G[1]$ vanishes in the 
Grothendieck group $\script K_0(V_k)$. 
Taking $G\set\R\mu_* \bigl(K_\bullet(\mathbf f_k)\bigr)^\sim$, we deduce then from Thomason's localization theorem\- \cite[p.\,338, 5.2.2(a)]{TT}   that the perfect $\SO_{V_k}$-complex 
$\R\mu_*C_k$ is $\mathbf D(V_k)$-isomorphic to the restriction of a perfect $\OX$-complex.
\begin{enumerate}
\item[(v)]
Let $S_k^{}\in\Dqc(X)$ be a perfect $\OX$-complex whose restriction
to $V_k$ is $\mathbf D(V_k)$-isomorphic to $\R\mu_*C_k$. \vs1
\item[(vi)]
Let $\cs_k$ be the finite set $\{S_1^{},S_2^{},\ldots,S_k^{}\}$.
\end{enumerate}

According to Lemma~\ref{Lemma3.3},  there is for each $k$ 
an integer $N_k>0$ such that, if $Q\in\D\qc(X)$ 
satisfies $H^\ell\<(Q)=0$ for all $\ell\ge -N_k$ then
\mbox{$\Hom_{\D(X)}(S_k^{},Q)=0$.} After enlarging $N_k$ if necessary, we have also 
that  $\Hom_{\D(X)}(Q,S_k^{})=0$.\vs{3.5}

Set

\begin{enumerate}
\item[(vii)]
$N\set\max\{\,N_1,N_2,\dots,N_t,r_1,r_2,\dots,r_t\,\}+1$.\vs1
\end{enumerate}
\entn

\medskip
Next comes the key statement.

\pro{P1.99}
With the preceding notation,  let\/ $m,k\in\mathbb Z,$ 
$1\le k\le t,$
let\/ $E\in\Dqc(X)$ be such that\/ $H^j(E)$ is supported in\/ $Y_k\>$ for all\/~$j\ge m-kN,$ and~set\vs{-2}
\[
a_k^{}=\binom{k+1}{2} N\qquad (1\le k\le t).
\]

{\rm (i)} If\/ $E$ is\/ $(m-(k-1)N)$-pseudo-coherent then there is
an\/ $m$-isomorphism\/ 
$P\to E$ with\/
$P\in\cs_k[m-a_k,\infty)$ \textup{(}so that\/ $P$ is perfect, see \textup{\ref{R1.0.3}(ii)).}\vs{1.5}

{\rm(ii)} If\/ $H^\ell\<(E)\ne0$ for some\/ $\ell\ge m,$ then for some\/ $i\ge m-a_k$ and 
some\/ $j\in[1,k],$ there is a nonzero map\/ $S_j[-i\>]\to E.$

\epro

Before proving this, let us see how to derive Theorems ~\ref{perfapprox} and ~\ref{perfgen}.

Since $Y_t=X\<$,  Theorem~\ref{perfapprox}, with $B\set (t-1)N$, is contained in 
~\ref{P1.99}(i).

Next,~\ref{P1.99}(ii) with $k=t$ shows that if $H^\ell\<(E)\neq0$, 
then there exist integers $i\geq \ell-a_t^{}$ and $j\in[1,t\>]$, 
and a non-zero map $S_j^{}[-i\>]\to E$. This gives Theorem~\ref{perfgen} for
the specific choices\vs{-3}
 $$
 S=S_1^{}\oplus S_2^{}\oplus\cdots\oplus S_t^{}\>,\qquad  A(S)\set  a_t=\binom{t+1}{2}N\,.
 $$

The rest of Theorem~\ref{perfgen} results from the following general fact, applied to  $\mathcal H=\{\,E\in\Dqc(X)\mid H^\ell\<(E)\ne 0\,\}$,  $\ct=\Dqc(X)$  and $A=A(S)-\ell$.

\pro{P1.17}
Let\/ $\ct$ be a triangulated category with coproducts. Let\/ $\mathcal H$~be a collection of objects of\/ $\ct$. Suppose there exists a compact generator\/
$S\in\ct$ and an integer $A$ such that\vs{-1}
$$
E\in\mathcal H\implies \Hom\bigl(S[n],E\bigr)\ne 0 \textup{ for some }n\le A.
$$
Then every compact generator has a similar property\/$:$ for each compact gener\-ator\/ $S'\in\ct$ there is an integer $A'$ such that\vs{-1}
$$
E\in\mathcal H\implies \Hom\bigl(S'[n\>],E\bigr)\ne 0 \textup{ for some }n\le A'.
$$
\epro

\prf
Let $\widehat R$ be the full subcategory of $\ct$ whose objects are all the direct summands of objects
in $\{S'\}(-\infty,\infty)=\bigcup_{M\ge0}\{S'\}[-M,M]$ (see Remark~\ref{R1.0.2}). 
As in the proof of Corollary~\ref{perf->pc}, (iii)${}\Rightarrow{}$(ii),  one sees that 
 $S\in\widehat R$, i.e.,
there is an $S^*\in\widehat R$ and
an~$M\ge 0$ such that $S\oplus S^*\in \{S'\}[-M, M]$.\vs1

Now if  $E\in\mathcal H$ then, since $\Hom\bigl(S[k],E\bigr)\ne 0$  for some $k\le A$,
and \mbox{$S[k]\oplus S^*[k]\in  \{S'\}[-M-k, \>M-k ]$} (Remark~\ref{R1.0.3}(i)), 
therefore\vs{.6}  Example~\ref{E1.0.5} gives
$\Hom(S'[n\>],E)\ne 0$ for some $n$ with $n\le M+k\le M+A$.\vs{-1}
\eprf

It remains to prove Proposition~\ref{P1.99}, which we do now by induction on $k$. \vs{1.5}

\enlargethispage*{2pt}

For $k=1$, $a_1=N\<$, so 
$H^j(E)$ is supported in $Y_1 = \Spec(R_1/I_1)$ for all~$j\ge m-r_1-1\;(\ge m-N)$. 
As usual, when considering the restriction $E|_{U_1}$ we may assume it to be
a quasi-coherent complex, then relate facts about it to facts about the corresponding complex~E of $R_1$-modules. For example, it holds that
$H^i$(E) is $I_1$-torsion for all~$i\ge m-r_1-1$.

Thus, from Corollary~\ref{L1.2.1}(i), 
applied to $I_1=(0,f_{11},f_{22},\cdots,f_{1r^{}_{\<\<1}})R_1$, it follows 
via~\ref{N1.3}(iv))
that, if $E$ is $m$-pseudo\kern.5pt-coherent then there exists an 
\mbox{$m$-isomorphism}
$P\to E|_{U_1}$ with $P\in\{C_1\}[m,\infty).$
Likewise (and more easily), 
Corollary~\ref{L1.2.1}(ii) gives that if 
$H^i\<(E)\ne 0$  for some $i\ge m$---whence, $H^i\<(E)$ being supported in $Y_1\subset U_1$, 
 $H^i\<(E|_{U_1)}\ne 0$---then there is a nonzero map\vs{-1}
$$
C_1[-i\>]= \big(K_\bullet(\mathbf f_1)[-i\>]\big)^\sim\oplus \big(K_\bullet(\mathbf f_1)[-i+1]\big)^\sim\to E|_{U_1}.
$$

Let $\mu\colon U_1\hookrightarrow X$ be the inclusion. Note that 
$Y_1=U_1\setminus V_2$. 
Since $C_1$ is exact outside $Y_1$, so is $P\in\{C_1\}[m,\infty)$ (argue as in Remark~\ref{R1.0.3}(i)--(iv)), as is $\R\mu_* P\>$; and by assumption, $H^i(E)$ vanishes outside~$Y_1$ for all  $i\ge m$. With all this in mind, 
we can extend the preceding statements from $U_1$ to ~$X=U_1\cup V_2$, by applying  the following Lemma to $U=U_1$, $V=V_2$, and $C=C_1$ or $P$.

\begin{lem}\label{extend} 
Let\/ $U$ and $V$ be open subsets of a scheme\/ $X\<,$ and let
$$
\CD
U\cap V @>\nu>> V \\
@V\lambda VV @VV\xi V \\
 U@>>\lift1.1,\mu,> U\cup V
 \endCD
 $$
be the natural diagram of inclusion maps. Let\/  $C\in\D(U)$ satisfy\/ 
$\lambda^*C=0$.  Let\/ $E\in\D(U\cup V)$. Then\/$:$\vs1 

\!{\rm (i)} Every\/ $\D(U)$-morphism\/ $C\<\to\mu^*E$
extends uniquely to a\/ $\D(U\cup V)$-morphism\/ \mbox{$\R\mu_*C\to E$.}\vs1

\!{\rm(ii)} If\/ $C$ is perfect\vs1 then so is\/ $\R\mu_*C$.

\!{\rm(iii)} If\/ $\cs\subset\D(U)$ and $m\le n\in\mathbb Z$ then 
$\R\mu_*\big(\cs[m,n]\big)\subset \big\{\R\mu_*\cs\big\}[m,n]$.

\end{lem}

\begin{proof}
(i) In view of the natural isomorphisms
$$
\Hom_{\D(U)}\<(C,\mu^*\<\<E)\<\cong\<\\
\Hom_{\D(U)}\<(\mu^*\R\mu_*C,\mu^*\<\<E)\<\cong\<
\Hom_{\D(U\cup V)}\<(\R\mu_*C,\R\mu_*\mu^*\<\<E)
$$
we need only show that the natural map is an isomorphism
$$
\R\Hom^\bullet(\R\mu_*C,E)\iso
\R\Hom^\bullet(\R\mu_*C,\R\mu_*\mu^*E)
$$
(to which we can apply the homology functor $H^0$). Thus for any triangle 
$$
\CD
G@>>> E @>\!\!\textup{natural}\!>>\R\mu_*\mu^*E@>>> G[1]
\endCD
$$
we'd like to see  that $\R\Hom^\bullet(\R\mu_*C,G)=0$. But  $\mu^*G=0=\lambda^*C$,
so that
\begin{align*}
&\mu^*\R\Hom^\bullet(\R\mu_*C,G)\cong
\R\Hom^\bullet(\mu^*\R\mu_*C,\>\mu^*G)=0,\ \textup{ and}\\
&\>\xi^*\>\R\Hom^\bullet(\R\mu_*C,G)\cong
\R\Hom^\bullet(\xi^*\R\mu_*C,\>\xi^*G)\cong
\R\Hom^\bullet(\R\nu_*\lambda^*C,\>\xi^*G)\!=\<0,
\end{align*}
whence the conclusion.\vs1

(ii) Since both $\xi^*\R\mu^*C\cong \R\nu_*\lambda^*C=0$ and $\mu^*\R\mu_*C=C$
are perfect, therefore so is $\R\mu_*C$.\vs1

(iii) This is a special case of Remark~\ref{R1.0.3}(iv).
\end{proof}

\begin{lem}\label{L1.12}
For\/ $k>1,$ suppose Proposition~\ref{P1.99}(i) holds with $k-1$ in place of\/  $k$. 
Then for any\/ $E\in\Dqc(X)$ 
and\/ $\D(U_k)$-morphism\vs{-2} 
\[
\psi\colon F\la E\vert_{U_k}^{}\qquad
\big(\> F\in\{C_k\}[m,\infty)\big),\vs{-2}
\]
there exists a\/ $\D(X)$-morphism\vs{-1}
\[
\tilde\psi\colon\wt F\to E \qquad
\big(\>\wt F\in \cs_k[m-N-a_{k-1}^{}\,,\,\infty)\>\big)\vs{-2}
\]
whose restriction\/
$\tilde\psi|_{U_k}$ is isomorphic to $\psi$.
\end{lem}

Before proving this Lemma, let us see how it is used to establish
the induction step in the proof of Proposition~\ref{P1.99}. With reference to that Proposition, we show, for $k>1$:

 (1) Assertion (i) for $k-1$ implies assertion (i) for $k$.

(2) Assertions (i) and (ii) for $k-1$, together, imply assertion (ii) for $k$.\vs2

To prove (1), let $E\in\Dqc(X)$ be $(m-(k-1)N)$-pseudocoherent, with
$H^j(E)$ supported in $Y_k$ for all $j\ge m-kN.$ Since
$m-(k-1)N-r_k\ge m-kN\<$, therefore (after replacement of 
$K_\bullet(\mathbf f_k)^\sim$ by $C_k$, see above) Corollary~\ref{L1.2.1}
provides a $\D(U_k)$-triangle
\[\label{tri}
\CD
P^{}_k @>>> E\vert_{U_k}^{} @>>> Q^{}_k @>>> P^{}_k[1]
\endCD\tag{\ref{L1.12}.1}
\]
with $P^{}_k\in \{C_k\}[m-(k-1)N\>\>,\,\infty)$ and $H^j(Q^{}_k)=0$ 
for all $j\geq m-(k-1)N$. 
By Lemma~\ref{L1.12}, the map $P^{}_k\la E\vert_{U_k}^{}$
is isomorphic\vs{.5}  to the restriction of a $\D(X)$-morphism $\psi'\colon P'\to E$, with
$P'\in \cs_k[m-(k-1)N-N-a_{k-1}^{}\,,\,\infty)$, i.e., since\looseness=-1 
\[
a_{k-1}^{}+kN\eq\binom{k}{2}N+kN\eq \binom{k+1}{2}N\eq
a_{k}^{}\>,\vs{-2}
\]
with $P'\in \cs_k[m-a_{k}^{},\infty)$. Any $\Dqc(X)$-triangle
\[\CD
P' @>\psi'>> E @>\alpha>> Q' @>>> P'[1]\,,
\endCD\]
restricts on $U_k$ to one isomorphic to \eqref{tri}.  So
when $j\geq m-(k-1)N$, then~$H^j(Q')$ vanishes on~$U_k$; furthermore,
$H^j(E)$ is supported on $Y_k\>$, and since all the members of $\>\cs_k$
are exact outside $Y_k$ therefore so is $P'$  
(argue as in~Remark~\ref{R1.0.3}(i)--(iv));
and thus $H^j(Q')$ is supported in $(Y_k\setminus
U_k)=Y_{k-1}\>$. 
\looseness=-1

Moreover, $Q'$ is $(m-(k-2)N\>)$-pseudocoherent, since both $P'$ and $E$
are \cite[p.\,100, 2.6]{I}.  
So now the inductive assumption produces
a  triangle
\[\CD
P'' @>>> Q' @>\beta>> Q @>>> P''[1]
\endCD\]
with $P''\in \cs_{k-1}[m-a_{k-1}^{}\>,\>\infty)$, and  $H^j(Q)=0$ whenever
$j\geq m$.

There is then a triangle\vs{-3}
\[\CD
P @>\psi'>> E @>\beta\alpha>> Q @>>> P[1]\,,
\endCD\]
and the assertion \ref{P1.99}(i), for the integer $k$,
results from Remark~\ref{R1.0.3}(v).\vs1

As for (2), let $E$ satisfy the hypotheses
of \ref{P1.99}(ii) for $k$.  If
$H^i(E\vert_{U_k}^{})=0$ for all 
$i\geq m-(k-1)N$ then   $H^j(E)$ is supported in $Y_{k-1}$
for all $j\geq m-(k-1)N$,  $H^\ell(E)$ is non-zero for some $\ell\geq m$,
and $m-a_{k-1}\ge m-a_k$; so in this case  assertion (ii) for $k$ is already
given by assertion (ii) for $k-1$.

If, on the other hand,\vs{.5}   
$H^i(E\vert_{U_k}^{})\neq 0$ for some $i\geq m-(k-1)N$, then, since 
$ m-(k-1)N-r_k\ge m-kN$, Corollary~\ref{L1.2.1}
(suitably modified) provides a
 nonzero map $C_k[-i\>]\to E\vert_{U_k}^{}$. 
By Remark~\ref{R1.0.3}(i), 
$$
C_k[-i\>]\in\{C_k\}[i,\infty)\subset \{C_k\}[m-(k-1)N,\infty),
$$
so by Lemma~\ref{L1.12}, there exists a nonzero $\D(X)$-morphism
$\wt F\to E$ with 
\[
\wt F\, \in\,\cs_{k}[m-(k-1)N-N-a_{k-1}^{}\,,\,\infty)
\;=\;
\cs_{k}[m-a_{k}^{},\infty) \: .
\]
Hence, by Example~\ref{E1.0.5}, \ref{P1.99}(ii) holds for $k$.

\medskip

 We come finally to the \emph{proof of Lemma~\ref{L1.12}.}\vs2

Let $\mathcal S\subset\D(U_k)$ be the full subcategory  with objects  those
$F\in\{C_k\}[m,\infty)$ for which the Lemma holds. We need to verify the conditions
in Definition~\ref{D1.0.1}, i.e., we need to show:\vs1

(a) $C_k[-\ell\>]\in\mathcal S$ for all $\ell\ge m$; and

(b) for any $\D(U_k)$-triangle\[
\CD
F' @>>> F @>>> F'' @>>> F'[1]\>,
\endCD
\]
if $F',F''\in\mathcal S$ then $F\in\mathcal S$.
\smallskip

For (a), we first use Lemma~\ref{extend} to extend\vs{-.5}
$\>\psi\colon C_k[-\ell\>]\to  E|_{U_k}$ to
a $\D(V_k)$-morphism $\phi\colon S_k[-\ell\>]\big|_{V_k}\to E|_{V_k}$.\vs1
By Thomason's localization
theorem, as formulated in \cite[p.\,214, 2.1.5]{N2} (and further
elucidated in [\emph{ibid.,} p.\,216, proof of Lemma 2.6]),%
\footnote{where in the absence of separatedness, $j_{\bullet*}$ should become 
$\R j_{\bullet*}$.}
 there is then
a $\Dqc(X)$-diagram, with top row  a triangle of perfect complexes:\vs{-1}\looseness=-1
\[
\CD
 \wt P @>>> \wt F_1 @>f>> S_k^{}[-\ell\>] @>>> \wt P[1]\\
@.  @V{g}VV @. @.\\
@.  E@. @.\quad.
\endCD
\]
and with $\wt P$ exact on $V_{k}$, so that $f|_{V_k}$ is an isomorphism; 
and furthermore,
$$
\phi=(g|_{V_k})\circ(f|_{V_k})^{-1}.
$$
Since $S_k[-\ell\>]\in\cs_k[\ell,\infty]$ (see Remark~\ref{R1.0.3}(i)), we need only show that we can choose 
$\wt P\in \cS_{k-1}[\ell-N-a_{k-1}^{}\,,\,\infty)$, because then we'll have\vs1
$$
\wt F_1 \in\cs_k[\ell-N-a_{k-1}^{}\,,\,\infty)\>\subset\>\cs_k[m-N-a_{k-1}^{}\,,\,\infty)\>.
$$

The perfect complex $\wt P\>$ is exact outside  $X-V_k=Y_{k-1}$, and we are 
assuming that \ref{P1.99}(i) is true for $k-1$. It follows that there exists a triangle
\[
\CD
 P @>>> \wt P @>>> Q @>>> P[1]
\endCD
\]
with $P\in \cS_{k-1}[\ell-N-a_{k-1}^{}\,,\,\infty)$ 
and $H^i(Q)=0$ for all $i\geq \ell-N$. Since all the members of $\cs_{k-1}$ are exact on $V_k$, the same is true of $P$  (argue as in Remark~\ref{R1.0.3}(i)--(iv)). 

\pagebreak[3]
Now \cite[p.\,58, 1.4.6]{N2} produces  
an octahedron on $P\to \wt P\to \wt F_1$, where the rows and columns are triangles:\vs2
\[
\CD
 P @= P \\
@VVV @VVV \\
\wt P @>>> \wt F_1 @>f>> S_k^{}[-\ell\>] @>>> \wt P[1]\\
@VVV @VV\beta V @| @VVV\\
Q @>>> F' @>>\lift1.3,f',> S_k^{}[-\ell\>] @>>\lift1,g,> Q[1]\\
@VVV @VVhV \\
 P[1] @= P[1] \vs2
\endCD
\]
Since $H^i\big(Q[1+\ell\>]\big)=0$ for all 
$i\geq  -N-1$,  the definition
of $N$ (see Notation~\ref{N1.3}(vii)) forces the
map $g\colon S_k^{}[-\ell\>] \to Q[1]$ to vanish. The exact sequence\vs{-2}
$$
\CD
\Hom(S_k[-\ell\>],\>F' )@>\textup{via}\;f'>> \Hom(S_k[-\ell\>],\>S_k^{}[-\ell\>] )
@>\!\textup{via}\;g\:=\:0\,>>
\Hom(S_k[-\ell\>], \>Q[1])\vs1
\endCD
$$
shows there is a map $\iota\colon S_k[-\ell\>]\to F'$ with
$f'\iota $  the identity map of $S_k^{}[-\ell\>]$.\vs1

This gives rise to yet another octahedron, on 
$S_k[-\ell\>]\xrightarrow{\iota} F'\xrightarrow{h}  P[1]$:\vs2
\[
\CD
 P @= P \\
@VVV @VVV \\ 
\wt F @>\gamma>> \wt F_1 @>>> Q @>>> \wt F[1]\\
@V{\alpha}VV @VV\beta V @| @VVV\\
S_k^{}[-\ell\>] @>\iota>> F' @>>> Q @>>> S_k^{}[-\ell+1]\\
@VVV @VVhV \\
 P[1] @= P[1] \vs{2.5}
\endCD
\]
The first column is a triangle,
with $\wt F \in\cs_k[\ell-N-a_{k-1}^{}\,,\,\infty)$,
and $P|_{V_k}$ exact, 
so that 
$\alpha|_{V_k}$~is an isomorphism. 

 Moreover, if $\>\tilde \psi\colon\wt F\to E$ 
is the composite $ \wt F\xrightarrow{\gamma} \wt F_1\xrightarrow{g} E$, then\vs{-1}
 $$
 \alpha=f'\iota\alpha=f'\beta\gamma=f\gamma,\vs{-3}
 $$ so that on $V_k$,\vs{-2}
 $$
 \phi\alpha=\phi f\gamma=g\gamma=\tilde\psi,\vs{-2}
 $$
proving (a).

\pagebreak[3]
\emph{Proof of} (b).

\medskip

Let $\psi\colon F\to E\vert_{U_k}^{}\<$ be a $\D(U_k)$-morphism.
Since $F'\in\mathcal S$, there exists a complex
$\>\wt F'\in \cS_k[m-N-a_{k-1}^{}\,,\,\infty)$ and
a $\D(X)$-morphism\vs{-4} $\wt{\psi'}\colon\wt F'\to E$\vs1 whose restriction to~$\>U_k$
is isomorphic to the composite\vs1 $F'\to F\xrightarrow{\psi}E\vert_{U_k}^{}$.
There results a triangle\vs{-2}\looseness=-1
\[
\CD
\wt F'@>\wt{\psi'}>>  E @>\gamma'>> E' @>>> \wt F'[1]\ ,\vs2
\endCD
\]
and hence a commutative $\D(U_k)$-diagram (part of an octahedron):\vs3
\[
\CD
F'@>>>  F @>>> F'' @>>> F'[1] \\
@|   @V\psi VV    @VV{g}V       @| \\
F' @>>> E\vert_{U_k}^{} @>\gamma'|_{U_k}>> E'\vert_{U_k}^{} @>>>F'[1] \\
@.      @V\chi VV    @VVhV @. \\
@.     G @= G\vs4   
\endCD
\]

Since $F''\in\mathcal S$, there is
an $\wt F''\in \cs_k[m-N-a_{k-1}^{}\,,\,\infty)$ and
a $\D(X)$-morphism $\wt{\psi''}\colon\wt F''\to E'$ whose restriction to $U_k$
is isomorphic to $g:F''\to E'\vert_{U_k}^{}$. So there is a triangle\vs{-3}
\[
\CD
\wt F''@>\wt{\psi''}>>  E' @>\gamma''>> E'' @>>> \wt F''[1]\ ;
\endCD
\]
whose restriction  to $U_k$ is isomorphic to\vs2
\[
\CD
F''@>g>>  E'\vert_{U_k}^{} @>h>> G @>>> F''[1]\ .\vs{-2}
\endCD
\]
The restriction to $U_k$ of the
composite\vs{.6} $E\xrightarrow{\gamma'} E'\xrightarrow{\gamma''} E''$ is isomorphic to
the composite $\chi\colon E\vert_{U_k}^{}\xrightarrow{\gamma'|_{U_k}} E'|_{U_k}\xrightarrow{h} G$.
Completing $\gamma''\gamma'$ to a triangle\vs{-3}
\[
\CD
\wt F@>\tilde\psi>>  E @>\gamma''\gamma'>> E'' @>>>\wt F[1]\vs1
\endCD
\]
and restricting to $U_k$,  we obtain  a triangle isomorphic to\vs1
\[
\CD
F@>\psi>>  E\vert_{U_k}^{} @>\chi>> G @>>>\wt F[1]\ .
\endCD
\]
That $\wt F\in \cs_k[m-N-a_{k-1}^{}\,,\,\infty)$
follows from Remark~\ref{R1.0.3}(v).
\qed

\pagebreak



%


\end{document}